\documentclass[a4paper, 12pt, draft, reqno]{amsart}
\usepackage[final]{graphicx}
\usepackage{amsmath}
\usepackage{amssymb}
\usepackage{latexsym}
\usepackage[cmtip,matrix,arrow]{xy}
\usepackage{rotating}
\UseComputerModernTips
\CompileMatrices

\newcommand{\ZZ}{{\mathbb Z}}

\newcommand{\RR}{{\mathbb R}}
\newcommand{\NN}{{\mathbb N}}
\newcommand{\CC}{{\mathbb C}}

\newcommand{\EE}{{\mathbb E}}
\newcommand{\HH}{{\mathbb H}}

\DeclareMathOperator{\Hom}{Hom}
\DeclareMathOperator{\Ext}{Ext}
\DeclareMathOperator{\End}{End}

\DeclareMathOperator{\res}{res}

\DeclareMathOperator{\supp}{supp}

\DeclareMathOperator{\id}{id}
\DeclareMathOperator{\spann}{span}

\DeclareMathOperator{\Alc}{\bf Alc}
\DeclareMathOperator{\MBS}{\bf MBS}

\DeclareMathOperator{\Reg}{Reg}
\DeclareMathOperator{\Graph}{\bf G}
\DeclareMathOperator{\Par}{\bf Par}
\DeclareMathOperator{\Mod}{-mod}

\newcommand{\otherwise}{\mbox{\rm otherwise}}
\newcommand{\wif}{\mbox{\rm if }}
\newcommand{\wand}{\mbox{\rm and }}

\DeclareMathOperator{\ind}{ind}

\DeclareMathOperator{\pr}{pr}

\newcommand{\too}{\longrightarrow}

\begin{document}
\theoremstyle{plain}
\newtheorem{thm}{Theorem}[section]
\newtheorem{prop}[thm]{Proposition}
\newtheorem{lem}[thm]{Lemma}
\newtheorem{cor}[thm]{Corollary}
\newtheorem{conj}[thm]{Conjecture}
\newtheorem{qn}[thm]{Question}
\newtheorem{claim}[thm]{Claim}
\theoremstyle{definition}
\newtheorem{rem}[thm]{Remark}
\newtheorem{ass}[thm]{Assumption}
\newtheorem{defn}[thm]{Definition}
\newtheorem{example}[thm]{Example}

\setlength{\parskip}{1ex}

\title[A translation principle for the Brauer algebra]{Alcove geometry
  and a  translation principle for the Brauer algebra} 
\author{Anton Cox} 
\email{A.G.Cox@city.ac.uk, M.Devisscher@city.ac.uk}
\author{Maud De Visscher} 
\address{Centre for Mathematical
  Science\\ City University\\ Northampton Square\\ London\\ EC1V
  0HB\\ England.}  \author{Paul Martin}
\email{ppmartin@maths.leeds.ac.uk} 
\address{Department of Pure
  Mathematics\\ University of Leeds\\ Leeds\\ LS2 9JT\\ England.}
 \begin{abstract}
There are similarities between algebraic Lie theory and a geometric
description of the blocks of the Brauer algebra in characteristic
zero. Motivated by this, we study the alcove geometry of a certain
reflection group action. We provide analogues of translation functors
for a tower of recollement, and use these to construct Morita
equivalences between blocks containing weights in the same
facet. Moreover, we show that the determination of decomposition
numbers for the Brauer algebra in characteristic zero can be reduced
to a study of the block containing the weight $0$. We define parabolic
Kazhdan-Lusztig polynomials for the Brauer algebra and show in certain
low rank examples that they determine standard module decomposition
numbers and filtrations.
 \end{abstract}

\maketitle

\section{Introduction and Motivation}

The Brauer algebra $B_n(\delta)$ was introduced in 1937 \cite{brauer}.
It is an enlargement of the symmetric group $\Sigma_n$, designed to
impose an orthosymplectic algebraic condition, depending on $\delta$,
on a Schur-Weyl dual algebraic group over the complex field.  However
it may be defined over an arbitrary ring $K$, for any $n \in \NN$ and
$\delta\in K$.  It has integral representations (in the sense of
\cite{bensonI}) that pass to simple modules over suitable splitting
fields (sometimes called cell modules in modern parlance), constructed
by Brown \cite{brownbrauer}.  This raises the problem of determining
simple decomposition matrices for these key modules (and hence for
indecomposable projective modules) over other extensions $k$ ($k$ a
field equipped with the property of $K$-algebra).  This long-standing
problem remains open.  The aim of the present work is twofold: to
develop some tools to solve this problem, by constructing a formal
`weight space' with a geometry and associated functors on the module
categories, and to propose a possible combinatorial framework (at
least over $\CC$) in which the answer might be couched.

The problem can be addressed in two parts: first working over $\CC$,
and then over fields of prime characteristic. (The latter can be
anticipated to be significantly harder, as the representation theory
of $B_n(\delta)$ {\em contains} the representation theory of the
symmetric group $\Sigma_n$.)  A significant step towards an answer
came with the determinantion in 2005 of the blocks of the algebra over
$\CC$ (and in 2006 of a geometric linkage principle in any
characteristic different from $2$) \cite{cdm,cdm2}.  These results
were obtained by using functors that allow the algebras for all $n$ to
be treated together (as previously used in \cite{mwdef,cgm}).  (An
alternative approach to the characteristic $0$ result via
characteristic $p$ has recently been developed by Donkin and Tange
\cite{dontan}.)

The key observation that underpins the various geometric
considerations in this paper is that the cell modules of $B_n(\delta)$
may be indexed by certain orbits of lattice points in the Euclidean
space $\EE^{\NN}$.  The orbits are those of a reflection group $A$,
where $A$ is the limit of the usual type-$A$ reflection group action
on Euclidean $N$-space $\EE^N$.  The $A$-action is a parabolic in the
limit type-$D$ reflection group action on $\EE^{\NN}$ (as in finite
rank), and the orbits of the $D$-action on coset space $\EE^{\NN}/A$
describe the blocks of $B_n(\delta)$-mod over $\CC$. For this reason
we will work over $\CC$ in this paper.  Using this
parabolic/reflection group formulation we are able to:
\begin{enumerate}
\item
determine a translation principle (Morita equivalences between certain
blocks);  
\item compute appropriate Brauer analogues of the
parabolic Kazhdan-Lusztig polynomials that determine decomposition
matrices in Lie (quantum group) theory; 
\item use these
to encode the structure of the algebra in many special cases (for
example in low rank, with the obvious conjecture that this extends to
all cases).
\end{enumerate}

Our methodology, and the structure of the paper, can be summarised as
follows.  The well-established root system/Weyl group analysis of
high-weight theory reduces many questions in Lie representation theory
(of algebraic and quantum groups) to geometry and combinatorics
\cite{deodhar,jannewed,soergel1}, once the Weyl group and affine Weyl
group action on weight space has been determined.  Of course these
Weyl groups are reflection groups, one a parabolic in the other
\cite{humcox}, facilitating, for example, an alcove geometric
description of blocks.

Note the obvious analogy with the role of reflection groups described
above.  It was this which motivated our formulation of the results in
\cite{cdm2} (guided by success with a similar approach to other
`diagram' algebras \cite{mwgen}).  In Lie theory the Euclidean
space is finite and the reflection group is infinite by virtue of
being affine; here it is by virtue of unbounded rank.  Nonetheless,
all the geometric and combinatorial machinery goes through unchanged.
The development of this analogy in Section \ref{weyl} lies at the
heart of our methodology.

Arguably one of the most beautiful machines that exists for computing
decomposition matrices in any setting is the method of (parabolic)
Kazhdan-Lusztig polynomials in Lie theory \cite{ajs,soergel1}.
Not all of the assumptions of this set-up hold for the Brauer algebra,
but in Section \ref{klp} we show how to bring the two theories close
enough together that parabolic Kazhdan-Lusztig (pKL) polynomials
suitable for the Brauer algebra may be computed.

In the group (or quantum group) case one has pKL polynomials
associated to alcoves in the alcove geometry, determining (at least
for the $q$-group over $\CC$, and its Ringel dual Hecke algebra
quotient \cite{erdqh}) decomposition matrices in alcove blocks. In
general there is more than one block intersecting an alcove, but there
is also a translation principle \cite{jannewed}, which states that all
these blocks are Morita equivalent, and hence do indeed have the same
decomposition matrices. In our case the pKL method formally assigns
the same decomposition matrix to every alcove block.  One is therefore
led to seek a form of translation principle.

The Brauer algebras $B_n(\delta)$ as $n$ varies form a tower via an
idempotent construction. In \cite{cmpx} we gave a general axiom scheme
for studying such a tower as a \emph{tower of recollement}. The
advantage of studying algebras in such a tower is the existence of
four functors: induction, restriction, globalisation, and
localisation, which relate the representation theories of the
different algebras in a compatible manner.

In Section \ref{towermor} we will show how towers of recollement, when
combined with a suitable description of the blocks in the tower, give
rise to analogues of translation functors and corresponding Morita
equivalences. These functors are defined using
induction or restriction functors followed by projection onto a block,
and are similar in spirit to $\alpha$-induction and
$\alpha$-restriction functors for the symmetric group \cite{gdeb}. 

We apply this translation theory to the Brauer algebras in Section
\ref{brauertran}, with the aim of proving that two blocks corresponding to
weights in the same facet have the same representation theory
(Corollary \ref{sumup}). However, in order to do this we will need
some additional functors, generalisations of induction and
restriction, which are introduced in Section \ref{genir}. We will also
see that when $\delta<0$ there are translation equivalences
\emph{between} certain facets, which raises interesting questions as
to the true geometric structure underlying the representation theory. 

We can also consider an analogue of translation \lq onto a wall\rq\ in
Lie theory for towers of recollement. Using this we show that the
decomposition matrix for the Brauer algebra is determined by the
decomposition matrix for the block containing the weight $0$ (Theorem
\ref{need0}).

In Section \ref{manygraphs} we consider various graphs associated to
each block (or each facet), and show that they are in fact all
isomorphic. For alcove graphs we can define associated Kazhdan-Lusztig
polynomials; using the graph isomorphisms these polynomials can more
generally be associated to any block graph.  In the final section we
show that when $\delta=1$ these polynomials correctly predict
decomposition numbers and filtrations in the alcove case for standard
modules in low rank examples.

\section{A review of Brauer algebra representation theory}
\label{basics}

In this section we will very briefly summarise the basic
representation theory of the Brauer algebras that will be needed in
what follows. Details can be found in \cite{cdm,cdm2}. In this paper,
we will restrict our attention to the case where the ground field is
$\CC$. 

The Brauer algebra $B_n(\delta)$ is a finite dimensional algebra with
parameter $\delta\in \CC$. When $\delta\notin \ZZ$ this algebra is
semisimple, so we will henceforth assume that $\delta$ is an
integer. \emph{We will also assume that $\delta\neq 0$.}

It will be convenient to use the usual graphical presentation of
Brauer algebras. An \emph{$(n,m)$ Brauer algebra diagram} will
consists of a rectangular frame with $n$ marked points on the northern
edge and $m$ on the southern edge called \emph{nodes}. Each of these
sets will be numbered from $1$ to $n$ (respectively $m$) from left to
right. Each node is joined to precisely one other by a line; lines
connecting the northern and southern edge will be called
\emph{propagating lines} and the remainder \emph{(northern or
  southern) arcs}.

Multiplication of two $(n,n)$ diagrams is by concatenation, where any
diagram obtained with a closed loop is set equal to $\delta$ times the
same diagram with the loop removed. Two diagrams are equivalent if
they connect the same pairs of nodes. The algebra obtained by taking
linear combinations of $(n,n)$ diagrams is a realisation of
$B_n$. Note that $\CC\Sigma_n$ is isomorphic to the subalgebra of $B_{n}$
spanned by diagram with only propagating lines. Moreover, $B_n$ is
generated by this subalgebra together with the elements $X_{i,j}$ with
$1\leq i,j\leq n$ consisting of $n-2$ propagating lines and arcs
joining $i$ and $j$ on the northern (respectively) southern edges.

The Brauer algebra can also be constructed via \lq iterated inflations' of
the symmetric group \cite{kxbrauer}, and thus is a cellular
algebra. If $\delta\neq 0$, then it is even quasihereditary. The standard
modules $\Delta_n(\lambda)$ are parameterised by partitions of
$n,n-2,\ldots,1/0$ (where the final term depends on the parity of
$n$), and we will denote the set of such by $\Lambda_n$. If
$\delta\neq 0$ then the same set parameterises the simple modules.

We have the following explicit construction of standard
modules. Consider a Brauer diagram with $n$ northern nodes and $n-2t$
southern nodes, and with no southern arcs. Such a diagram must have
exactly $t$ northern arcs. We will denote this diagram by
$X_{v,1,\sigma}$, where $v$ denotes the configuration of northern
arcs, $1$ represents the fixed southern boundary, and
$\sigma\in\Sigma_{n-2t}$ is the permutation obtained by setting
$\sigma(i)=j$ if the $i$th propagating northern node from the left is
connected to the southern node labelled by $j$.

 The elements $v$ arising as above will be called \emph{partial
   one-row diagrams}, and the set of such will be denoted by
 $V_{n,t}$. If a node $i$ in $w\in V_{n,t}$ is not part of a northern
 arc we say that it is \emph{free}.  The vector space spanned by the
 set of diagrams of the form $X_{w,1,id}$ where $w\in V_{n,t}$ will be
 denoted $I_n^t$. Note that $\Sigma_{m}$ acts on $I_n^t$ on the right
 by permuting the southern nodes.

Given $\lambda$ a partition of $m=n-2t$, let $S^{\lambda}$ denote the
Specht module corresponding to $\lambda$ for $\Sigma_{m}$. Then the
standard module $\Delta_n(\lambda)$ can be realised (see \cite[Section
  2]{dhw} or \cite[Section 2]{cdm}) in the following manner. As a
vector space we have
\begin{equation}\label{concrete}
\Delta_{n}(\lambda)=I_n^t\otimes S^{\lambda}.
\end{equation}
An element $b$ of $B_n$ acts on $d\in I_n^t$ from the right by diagram
multiplication. If the resulting product has fewer than $m$
propagating lines then we define the action of $b$ on $d\otimes
S^{\lambda}$ to be $0$. Otherwise the
product will result in a diagram with exactly $m$ propagating lines,
but these may now be permuted. We transfer this permutation (thought
of as an element of $\Sigma_m$) through the tensor product to act on
$S^{\lambda}$.

For $\delta\neq 0$ (or $n>2)$ there is an idempotent $e_n\in
B_{n}(\delta)$ such that $e_nB_n(\delta)e_n\cong B_{n-2}(\delta)$, and
so there are associated localisation and globalisation functors
$F_n:B_{n}(\delta)\mbox{\rm -mod} \rightarrow B_{n-2}(\delta)\mbox{\rm
  -mod}$ and $G_n:B_{n}(\delta)\mbox{\rm -mod} \rightarrow
B_{n+2}(\delta)\mbox{\rm -mod}$. In this way we can regard
$B_{n}(\delta)\mbox{\rm -mod}$ as a full subcategory of
$B_{n+2}(\delta)\mbox{\rm -mod}$, and hence
$\Lambda_n\subset\Lambda_{n+2}$. We set
$\Lambda=\lim_{n\rightarrow\infty}(\Lambda_n\cup\Lambda_{n+1})$, the
set of all partitions.  We will abuse terminology and say that two
labels are in the same block when the associated standard modules are in
the same block.

In order to describe the main results in \cite{cdm} we will need some
additional terminology. Recall that for a partition $\lambda$ (which
we will identify with its Young diagram), the \emph{content} of the
box in row $i$ and column $j$ of the diagram is defined to be
$j-i$. 
A pair of partitions $\mu\subset\lambda$ is said to be
\emph{$\delta$-balanced} if the following conditions are satisfied:
\begin{enumerate}
\item  the boxes in the skew partition can be paired so
that the contents of each pair sum to $1-\delta$;
\item if the skew partition contains boxes labelled by 
$1-\frac{\delta}{2},-\frac{\delta}{2}$ and there is only one such box
  in the bottom row then the number of pairs of such boxes is even.
\end{enumerate}
We say that two general partitions $\lambda$ and $\mu$ are
\emph{$\delta$-balanced} if the pairs $\lambda\cap\mu\subset \lambda$ and
$\lambda\cap\mu\subset\mu$ are both $\delta$-balanced.
The importance of the $\delta$-balanced condition is clear from the 
 following result \cite[Corollary 6.7]{cdm}:

\begin{thm}\label{balver}
Two partitions $\lambda$ and $\mu$ are in the same block for
$B_n(\delta)$ if and only if they are $\delta$-balanced.
\end{thm}

Denote by $V_{\delta}(\lambda)$ the set of partitions $\mu$ such that
$\mu$ and $\lambda$ are $\delta$-balanced. Note that if $\mu\in
V_{\delta}(\lambda)$ then so too are $\lambda\cap\mu$ and
$\lambda\cup\mu$. Thus $V_{\delta}(\lambda)$ forms a lattice under the
inclusion relation. We say that $\mu$ is a \emph{maximal balanced
  subpartition of} $\lambda$ if $\mu\in V_{\delta}(\lambda)$ and there
does not exist $\tau\in V_{\delta}(\lambda)$ with
$\mu\subset\tau\subset\lambda$.  One of the main steps in the proof of
Theorem \ref{balver} is \cite[Theorem 6.5]{cdm}, which shows that if
$\mu$ is a maximal balanced subpartition of $\lambda$ then
$$\Hom_n(\Delta_n(\lambda),\Delta_n(\mu))\neq 0.$$

The standard module $\Delta_n(\lambda)$ has simple head $L_n(\lambda)$
and all other composition factors are of the form $L_n(\mu)$ where
$\mu\supset\lambda$ lies in the same block as $\lambda$
\cite[Proposition 4.5]{cdm}. If $\lambda$ and $\mu$ are such a pair
with $|\lambda/\mu|=2$ then
\begin{equation}\label{twoboxstuff}
[\Delta_n(\mu):L_n(\lambda)]=\dim\Hom(\Delta_n(\lambda),\Delta_n(\mu))=1
\end{equation}
by  \cite[Theorem 3.4 and the remarks after Theorem 3.1]{dhw}
(see \cite[Theorem 4.4]{cdm}). If $\mu\subset \lambda$ are two weights
in the same block and $\lambda/\mu=(a^b)$ for some $a$ and $b$ with
$a$ even then we
also have \cite[Proposition 5.1]{cdm} that
\begin{equation}\label{squarestuff}
[\Delta_n(\mu):L_n(\lambda)]=1.
\end{equation}
In general, if $\lambda\vdash n$ and $\mu\vdash m$ with $m\leq n$ the
exactness of the localisation functor implies that
$$[\Delta_N(\mu):L_N(\lambda)]=[\Delta_n(\mu):L_n(\lambda)]$$
for all $N>n$.

The algebra $B_n(\delta)$ embeds inside $B_{n+1}(\delta)$, and so we
may consider the associated induction and restriction functors
$\ind_n$ and $\res_{n+1}$. If $\lambda$ and $\mu$ are partitions we
write $\lambda\rhd\mu$, or $\mu\lhd\lambda$ if the Young diagram for
$\mu$ is obtained from that for $\lambda$ by removing one box. Then
\cite[Theorem 4.1 and Corollary 6.4]{dhw} we
have short exact sequences
\begin{equation}\label{indrule}
0\rightarrow \bigoplus_{\mu \lhd \lambda} \Delta_{n+1}(\mu)
\rightarrow {\rm ind}_n\, \Delta_n(\lambda) \rightarrow \bigoplus_{\mu
  \rhd \lambda} \Delta_{n+1}(\mu)\rightarrow 0
\end{equation}
and
\begin{equation}\label{resrule}
0\rightarrow \bigoplus_{\mu \lhd \lambda} \Delta_{n-1}(\mu)
\rightarrow {\rm res}_n\, \Delta_n(\lambda) \rightarrow \bigoplus_{\mu
  \rhd \lambda} \Delta_{n-1}(\mu)\rightarrow 0.
\end{equation}

The restriction rule for simples is not so straightforward. However we
do have by \cite[Lemma 7.1]{cdm} that if $\mu$ is a partition obtained
from $\lambda$ by removing one box then 
\begin{equation}
\label{ressimple}
[\res_n L_n(\lambda):L_{n-1}(\mu)]\neq 0.
\end{equation}

The next two results are new, and show how the local data in
(\ref{indrule}) and (\ref{resrule}) can be applied to explicit
decomposition number calculations (which illustrates one of the motivations
for the tower of recollement formalism in \cite{cmpx}).

\begin{prop}\label{projelim}
Suppose that by removing $m$ boxes from $\lambda\vdash n$ it is
possible to reach a partition $\mu\vdash n-m$ such that $\Delta_{n-m}
(\mu)$ is a projective $B_{n-m}(\delta)$-module. Then the simple
module $L_n(\lambda)$ does not appear as a composition factor in any
$\Delta_{n}(\nu)$ with $\nu$ of degree less than $n-2m$.
\end{prop}
\begin{proof} Suppose that $L_n(\lambda)$ does occur as a composition
  factor of $\Delta_n(\nu)$. Then by Brauer-Humphreys reciprocity
  \cite[Proposition A2.2(iv)]{don2} the projective cover
  $P_n(\lambda)$ of $L_n(\lambda)$ has a standard module filtration
  with $\Delta_n(\nu)$ as a factor.

As $\Delta_{n-m}(\mu)$ is projective, so is
$\ind_{n-1}\cdots\ind_{n-m}\Delta_{n-m}(\mu)$. By repeated application
of (\ref{indrule}) we see that this contains $L_n(\lambda)$ in its
head, and so must have as a summand $P_n(\lambda)$. We have that
$$\ind_{B_{n-m}}^{B_n}=\res_{B_n}^{B_{m+n}}G_{n+m-2}\cdots
G_{n-m+2}G_{n-m}$$
and 
$$G_r(\Delta_r(\lambda))\cong \Delta_{r+2}(\lambda).$$
Therefore by repeated application of (\ref{resrule}) to
$\Delta_{n+m}(\mu)$ we see that $P_n(\lambda)$
cannot have a standard module filtration with $\Delta_n(\nu)$
as a factor. This gives the desired contradiction and so we are done.
\end{proof}

\begin{rem}
Note that any standard module which is alone in its block must be
projective. Thus there are many circumstances where Proposition
\ref{projelim} will be easy to apply. Indeed, this case will be
sufficient for our purposes.
\end{rem}

If $\mu\subset\lambda$ are two partitions then their skew
$\lambda/\mu$ can be regarded as a series of disjoint
partitions; when considering such differences we will list the various
partitions in order from top right to bottom left. Thus a skew
partition $((22)^2)$ will consist of two disjoint partitions of the
form $(22)$. 

\begin{prop}\label{largerhom}
If $\mu\subset\lambda$ is a balanced pair with $\lambda/\mu=((22)^2)$
or $\lambda/\mu=((1)^4)$
then 
$$\Hom_n(\Delta_n(\lambda),\Delta_n(\mu))\neq 0.$$
\end{prop}
\begin{proof}
We may assume that $\lambda\vdash n$ by localisation.  If
$\lambda/\mu=((22)^2)$ let $\lambda'$ be $\lambda$ less one of the two
removable boxes in $\lambda/\mu$, and $\mu'$ be the partition $\mu$
together with the addable box from the other component of the
skew. If $\lambda/\mu=((1)^4)$ then let $\lambda'$ be $\lambda$ with
any one of the boxes in $\lambda/\mu$ removed, and $\mu'$ be $\mu$
together with the unique box in $\lambda/\mu$ making this a balanced pair.

In each case  $\mu'$ is a maximal balanced subpartition of
$\lambda'$ and so by \cite[Theorem 6.5]{cdm} we have that
$$\Hom_{n-1}(\Delta_{n-1}(\lambda'),\Delta_{n-1}(\mu'))\neq 0.$$ By
\cite[Corollary 6.7]{cdm} and (\ref{resrule}) the only term in the
block labelled by $\lambda'$ in the standard filtration of
$\res_{n}\Delta_{n}(\mu)$ is $\Delta_n(\mu')$. Therefore by
Frobenius reciprocity we have
\begin{equation}\label{frex}
\Hom(\ind_{n-1}\Delta_{n-1}(\lambda'),\Delta_n(\mu))\cong
\Hom(\Delta_{n-1}(\lambda'),\res_n\Delta_n(\mu))\neq 0.
\end{equation}
Now $\lambda$ is not the only weight in its block in the set of weights
labelling term in the standard filtration of
$\ind_{n-1}\Delta_{n-1}(\lambda')$. However, by \cite[Lemma 4.10]{cdm}
it follows from (\ref{frex}) that 
$$\Hom_n(\Delta_n(\lambda),\Delta_n(\mu))\neq 0$$
as required.
\end{proof}

In \cite{cdm2} we identified partitions labelling Brauer algebra
modules with elements of $\ZZ^N$ (for suitable $N$) only after
transposition of the original partition to form its
conjugate. Henceforth when we regard $\Lambda_n$ (or
$\Lambda$) as a subset of $\ZZ^{\infty}$ it will always be via this
\emph{transpose map} $\lambda\rightarrow \lambda^T$.

\section{Brauer analogues of Weyl and affine Weyl groups}\label{weyl}

We wish to identify reflection groups associated to the Brauer algebra
which play the role of the Weyl and affine Weyl groups for reductive
algebraic groups. First let us recall the properties of Weyl groups
which we wish to replicate.

In Lie theory a Weyl group $W$ is a reflection group acting on a
Euclidean weight space with the following properties:
\begin{enumerate}
\item There is an integral set of weights on which $W$ acts via a
  \lq dot' action $(W,\cdot)$,
\item The reflection hyperplanes of $W$ under this action break space
  up into chambers (and other facets),
\item A complete set of weights indexing simple (or standard) modules
  coincides with the weights in a single chamber under the dot action,
  namely that containing the zero weight. Such weights are said to be
  \emph{dominant}.
\end{enumerate}

Thus the selection of an indexing set for the dominant weights is
taken care of by the Weyl group (and its dot action). In positive
characteristic $p$ or at a quantum $l$th root of unity there is then a
second stage, the introduction of an affine extension of $W$ (with
action depending on $p$ or $l$), which has orbits whose intersection
with the dominant weights determine the blocks.

This affine extension defines an additional set of reflecting
hyperplanes, which break the set of weights up into a series of
chambers (now called alcoves) and other facets. We refer to this
configuration of facets, together with the action of the affine
extension, as the \emph{alcove geometry} associated to the particular
Lie theory in question.  

The alcove geometry controls much of the representation theory of the
corresponding reductive group. In particular, we typically have a
\emph{translation principle} which says that there are Morita
equivalences between blocks which intersect a given facet, and so much
of the representation theory does not depend on the weight itself but
only on the facet in which it lies.

We will show how a version of the above programme can be implemented 
for the Brauer algebra from scratch.

Let $\EE^n$ be the $\RR$-vector space with basis $e_1,\ldots, e_n$,
and $\underbar{n}=\{1,\ldots, n\}$. We will define various reflections
on $\EE^n$ corresponding to the standard action of the type $D$ Weyl
group. Let $(ij)$ be the reflection in the hyperplane in $\EE^n$
through the origin which takes $e_i$ to $e_j$ and fixes all other unit
vectors, and $(ij)_-$ to be the reflection in the hyperplane
perpendicular to $e_i+e_j$ which takes $e_i$ to $-e_{j}$.  We define
$$W_a(n)=\langle(i,j),(i,j)_- : i\neq j\in\underbar{n}\rangle$$
which is the type $D$ Weyl group. Note that it has a 
a subgroup
$$W(n)=\langle(i,j) : i\neq j\in\underbar{n}\rangle$$
which is just the type $A$ Weyl group (isomorphic to $\Sigma_n$).

As explained in the previous section, the Brauer algebras
$B_n(\delta)$ as $n$ varies form a tower. Thus it is natural to
consider all such algebras simultaneously. In order to do this we will
work with the infinite rank case. Note also that orbits of the finite
Weyl group $W(n)$ are not sufficient to define an indexing set for the
simple $B_n(\delta)$-modules (one needs to consider $W(n+1)$-orbits,
but then this group is not a subgroup of $W_a(n)$), unlike the
infinite rank case.

Let $\EE^{\infty}$ be the
$\RR$-vector space consisting of (possibly infinite) linear
combinations of the elements $e_1,e_2,\ldots$. We say that
$\lambda\in\ZZ^{\infty}$ has \emph{finite support} if only finitely
many components of $\lambda$ are non-zero, and write $\ZZ^f$ for the
set of such elements.  (We define $\EE^f$ similarly.) Thus, with the
obvious embedding of $\ZZ^n$ inside $\ZZ^{n+1}$, we have that
$\ZZ^f=\lim_{n\rightarrow\infty}\ZZ_n$.  (Note that the usual
component-wise inner product on $\ZZ^{\infty}$ is finite on $\ZZ^f$, so
that $\ZZ^f$ is just the integral part of the usual $l^2$ Hilbert
space in $\RR^{\infty}$.)

We say that an element $\lambda=(\lambda_1, \lambda_2,\ldots)$ of
$\ZZ^{f}$ is \emph{dominant} if $\lambda_i\geq\lambda_{i+1}$ for all
$i$. (Note that any such element must lie in $\NN^{\infty}$.) Embed
$\EE^n$ inside $\EE^{\infty}$ in the obvious way, and let $W_a$, and
$W$ be the corresponding limits of $W_a(n)$, and $W(n)$. Clearly the
space $\ZZ^f$ is closed under the action of $W_a$. We will call
elements in $\ZZ^f$
\emph{weights}. Dominant weights are precisely those which
label standard modules for the Brauer algebras, and by analogy with Lie
theory we will denote the set of such weights by $X^+$. General
elements of $\EE^{\infty}$ will be called \emph{vectors}. We will use
Greek letters for weights and Roman letters for general vectors.

Given a reflection group $G$ (or the corresponding set of hyperplanes
$\HH$), we say that a vector is \emph{regular in $G$ (or in $\HH$)} if
it lies in the interior of a chamber, i.e. in some connected component
of $\EE^{\infty}\backslash \cup_{X\in\HH}X$. Otherwise we say the
vector is \emph{singular}. In the case $G=W_a$ we shall call chambers
\emph{alcoves} to emphasise the distinction between this and the
$W$ case. For $v\in\EE^{\infty}$ we define the \emph{degree of
  singularity}
$$s(v)=|\{\{i,j\}:v_i=\pm v_j, \quad i\neq j\}|$$
(which need not be finite in general).  
Note that a vector $v$ is regular in $W_a$ if and only if
$s(v)=0$. The next lemma is clear.

\begin{lem}\label{basicalcoves}
(i) There is a chamber $A^+$ of the action of $W$ on
  $\EE^{\infty}$ consisting of all strictly decreasing
  sequences.\\ 
(ii) The boundary of $A^+$ consists of all non-strictly
  decreasing sequences.
\end{lem}

Recall that in Lie theory we typically consider a shifted reflection
group action with respect to some fixed element $\rho$. It will be
convenient to consider a similar adjustment here. Let
$-2\omega=(1,1,\ldots,1)\in\EE^{\infty}$ and
$\rho_0=(0,-1,-2,\ldots)$.
For $\delta\in\ZZ$ define
$$\rho_\delta=\rho_0+\delta\omega.$$

For $w\in W_a$ and $v\in\EE^{\infty}$ let
$$w\cdot_\delta v=w(v+\rho_{\delta})-\rho_{\delta}$$ where the
right-hand side is given by the usual reflection action of $W_a$ on
$\EE^{\infty}$. Note that $\ZZ^f$ is closed under this action of
$W_a$.  We say that a weight $\lambda$ is \emph{$\delta$-regular} if
$\lambda+\rho_{\delta}$ is regular, and define the \emph{degree of
  $\delta$-singularity} of $\lambda$ to be $s(\lambda+\rho_{\delta})$.

\begin{prop}\label{Wis} Let $\lambda\in\ZZ^f$.\\
(i) For $w\in W$ the weight
  $w\cdot_{\delta}\lambda$ does not depend on $\delta$. Moreover, if
  $w\neq 1$ and $\lambda\in X^+$ then $w\cdot_{\delta}\lambda\notin X^+$.\\
(ii) If $\lambda\in X^+$ then $\lambda+\rho_{\delta}$ can only lie on
  a $(ij)_-$-hyperplane.\\
(iii) We have $\lambda\in X^+$ if and only if
  $\lambda+\rho_{\delta}\in A^+$.
\end{prop}
\begin{proof}
(i) Note that 
$$(ij)(\lambda+\rho_{0}+\delta\omega)-\rho_{0}-\delta\omega
=(ij)(\lambda+\rho_0)-\rho_{0}.$$ 
(ii) and (iii) are clear.
\end{proof}

The description of the blocks of the Brauer algebra in characteristic
zero in Theorem \ref{balver}  was given the following geometric
reformulation in \cite{cdm2}:

\begin{thm}\label{geoblock}
Two standard modules $\Delta_n(\lambda)$ and $\Delta_n(\lambda')$ for
$B_n(\delta)$  are in the same block if and only if
$\lambda^T$ and $\lambda'^T$ are in the same
$(W_a(n),\cdot_{\delta})$-orbit.
\end{thm}

\begin{rem}
In summary, we have shown that there is a space that plays a role
analogous to a weight space (in Lie theory) for the Brauer algebra,
together with an action of the type $A$ Coxeter group which plays the
role of the Weyl group, while the corresponding type $D$ Coxeter group
plays the role of the affine Weyl group.
\end{rem}

We will now consider the geometry of facets induced by $W_a$ inside
$A^+$. We will call reflection hyperplanes \emph{walls}, and for any
collection of hyperplanes we will call a connected component of the
set of points lying on the intersection of these hyperplanes but on no
other a \emph{facet}. (Then an alcove is a facet corresponding to the
empty collection of hyperplanes.)

It will be convenient to have an explicit description of the set of
vectors in a given facet.
For vectors in $A^+$ (which will be the only ones which concern us)
these facets are defined by the hyperplanes $v_i=- v_j$ for some
$i\neq j$. For $v=(v_1, v_2, v_3, ...)$ in $A^+$, note that for all
$i\in\NN$ we have
$$|\{j \, : \, |v_j|=|v_i|\}| \leq 2.$$ We will call $v_i$ a
\emph{singleton} if $v_i$ is the only coordinate with modulus $|v_i|$,
and the pair $v_i, v_j$ a \emph{doubleton} if $|v_i|=|v_j|$ and $i\neq
j$.

For a given facet $F$ with $v\in F$, a vector $v'\in A^+$ lies in $F$
if and only if $|v'_i|=|v'_j|$ whenever $|v_i|=|v_j|$, and
$|v'_i|>|v'_j|$ whenever $|v_i|>|v_j|$. Therefore an alcove (where
every $v_i$ is a singleton) is determined by a permutation $\pi$ from
$\NN$ to $\NN$ where $|v_{\pi(n)}|$ is the $n$th smallest modulus
occurring in $v$. Note that not every permutation corresponds to an
alcove in this way.  Further, if $i<\pi(1)$ then $v_i>0$,
while if $i>\pi(1)$ then $v_i<0$.

For more general facets we replace the permutation $\pi$ by a function
$f:\NN\rightarrow \NN\cup (\NN\times\NN)$ such that $f(n)$ is the
coordinate, or pair of coordinates, where the $n$th smallest modulus in
$v$ occurs. For example, if 
$$v=(6,4,2,1,0,-2,-3,-5,\ldots)$$
then the facet containing $v$ corresponds to a function whose first
four values are $f(1)=5$, $f(2)=4$, $f(3)=(3,6)$, and $f(4)=7$.

We will denote by $A_0$ the alcove corresponding to the identity
permutation. Thus $A_0$ consists of all $v\in A^+$ such that
$|v_1|<|v_2|$ and $v_2<0$. It is easy to see that, for any $\delta\geq
0$, the weight $0$ is $\delta$-regular, with the vector
$0+\rho_{\delta}$ in $A_0$; in this case we will call the
$(W_a,\cdot_{\delta})$-alcove the \emph{$\delta$-fundamental alcove}.

\begin{lem}\label{poscase}
For $\delta\geq 0$ the set of weights in the $\delta$-fundamental
alcove is 
$$\{\lambda\in X^+:\lambda_1+\lambda_2\leq \delta\}.$$
\end{lem}
\begin{proof}
By our discussion above, the desired set of weights is precisely the
set of dominant $\lambda$ such that $x=\lambda+\rho_{\delta}$
and $|x_1|<|x_2|$. But this means that
$$\lambda_1-\frac{\delta}{2}<\frac{\delta}{2}+1-\lambda_2$$
which implies the result.
\end{proof}

\begin{rem}
Although our alcove geometry is reminiscent of that arising in
positive characteristic Lie theory, there are also some striking
differences. Consider for example the case when $\delta=1$. The alcove
$A_0$ is non-empty and contains the two weights $0$ and $(1)$. The
next lowest alcove contains $(2,1)$ and $(2,2)$, and the third
contains $(3,2,1)$ and $(3,1,1)$. However, the associated facets with
singularity $1$ are not necessarily finite (in $\EE^{\infty}$); for
example there is a facet consisting of the weights $(2)$ and the
$n$-tuple $(1,\ldots, 1)$ for all $n\geq 2$. In particular, not every
weight on a wall is adjacent to a weight in an alcove.
\end{rem}

Recall that for $\delta\in \NN$ the Brauer algebra $B_r(\delta)$ is in
Schur-Weyl duality with O$_{\delta}(\CC)$ acting on the $r$th tensor
product of the natural representation. 

\begin{thm}
Suppose that $\delta\in\NN$. The elements of $\Lambda_{\infty}$
corresponding to weights in the $\delta$-fundamental alcove are in
bijection with the set of partitions labelling the irreducible
representations which arise in a decomposition of tensor powers of the
natural representation of \mbox{\rm O}$_n(\CC)$.
\end{thm}
\begin{proof}
For O$_n(\CC)$ tensor space components are labelled by partitions
whose first and second columns sum to at most $n$ (see for example
\cite[Theorem 10.2.5]{goodwall}). The result now follows by
comparing with Lemma \ref{poscase} via the transpose map on partitions.
\end{proof}

\begin{rem}
The above result shows that the fundamental alcove arises naturally in
the representation theory of O$_n(\CC)$.
\end{rem}

Suppose that $\delta<0$. Choose $m\in\NN$ so that $\delta=-2m$ (if
$\delta$ is even) or $\delta=-2m+1$ (if $\delta$ is odd). It is easy
to see that $0$ is $\delta$-singular of degree $m$. Indeed, any
dominant weight $\lambda$ is $\delta$-singular of degree at least
$m$. Thus there are no regular dominant weights for
$\delta<0$. Instead of the $\delta$-fundamental alcove, we can
consider the \emph{$\delta$-fundamental facet} containing $0$, for
which we have

\begin{lem}\label{negcase}
For $\delta< 0$ of the form $-2m$ or $-2m+1$ the set of weights in the
$\delta$-fundamental facet is $\{0\}$.
\end{lem}
\begin{proof}
We consider the case $\delta=-2m$; the odd case is similar. The
element $0+\rho_{\delta}$ equals
$$(m,m-1,\ldots,0,-1,\ldots)$$ 
and hence our facet consists of all vectors of the form
$$(t,t-1,\ldots,-t+1,-t,v_{|\delta|},v_{|\delta|+1},\ldots)$$ 
where the sequence $-t,v_{|\delta|}\ldots$ is decreasing (as any
other weight would be non-dominant).
But this implies that $Y_{\delta}=\{0\}$.
\end{proof}

This is very different from the case $\delta>0$.
However we do have

\begin{thm}\label{fakealc}
 Suppose that $\delta=-2m$. The set of elements of
  $\Lambda_{\infty}$ corresponding to weights which can be obtained
  from $0$ via a sequence of one box additions only involving
  intermediate weights
  of singularity $m$ is in bijection with the set of partitions
  labelling the irreducible representations which arise in a
  decomposition of tensor powers of the natural representation of
  \mbox{\rm Sp}$_{2m}(\CC)$.
\end{thm}
\begin{proof} First note that we can clearly add boxes in the first
  $m$ coordinate of $0$ without changing the degree of
  singularity. In order to change the $(m+1)$st coordinate in our path,
  we will have to pass through some point of the form
$$(a_1,a_2,\ldots,a_m,1,-1,-2, \ldots).$$ where $a_1> a_2>\cdots>a_m>1$. But
this implies that the first $m+1$ coordinates of the vector all pair up
with the corresponding negative values later down the vector, and so
this is a singular vector of degree $m+1$. The result now follows from
the description of tensor space components (see for example
\cite[Theorem 10.2.5]{goodwall}).
\end{proof}

We will see in Section \ref{brauertran} that there is a sense in which
the set of weights occurring in Theorem \ref{fakealc} can be regarded
as playing the role of an alcove in the $\delta<0$ case.

\section{A translation principle for towers of recollement}\label{towermor}

Towers of recollement were introduced in \cite{cmpx} as an axiom
scheme for studying various families of algebras. The Brauer algebra
over $\CC$ was shown to satisfy these axioms in \cite{cdm}. We will
prove a general result about Morita equivalences in such towers, and
apply it in the following section to the Brauer algebra. In this
section we will work over a general field $k$.

Suppose that we have a family of $k$-algebras $A=\{A_n:n\in\NN\}$
forming a tower of recollement. The precise properties of such a tower
will not concern us here; such properties as we need will be
introduced in what follows, and details can be found in \cite{cmpx}.

In each $A_n$ with $n\geq 2$ there exists an idempotent $e_n$ such
that $e_nA_ne_n\cong A_{n-2}$. This determines a pair of functors:
\emph{localisation} $F_n$ from $A_n\Mod$ to $A_{n-2}\Mod$, and
\emph{globalisation} $G_n$ from $A_n\Mod$ to $A_{n+2}\Mod$ given on
objects by
$$F_nM=e_nM\quad\quad\wand\quad\quad
G_nM=A_{n+2}e_{n+2}\otimes_{A_n}M.$$ 
The functor $F_n$ is exact, $G_n$
is right exact, and $G_n$ is left adjoint to $F_{n+2}$. 
We also have algebra inclusions $A_n\subset A_{n+1}$ for
each $n\geq 0$, and associated functors $\ind_n$ from
$A_n\Mod$ to $A_{n+1}\Mod$ and $\res_n$ from $A_n\Mod$ to $A_{n-1}\Mod$.
Let $\Lambda_n$ be an indexing set for the simple
$A_n$-modules. Globalisation induces an embedding of $\Lambda_n$
inside $\Lambda_{n+2}$, and we take $\Lambda$ to be the disjoint union
of $\lim_{n}\Lambda_{2n}$ and $\lim_n\Lambda_{2n+1}$, whose elements
we call \emph{weights}.  

The algebras $A_n$ are quasihereditary, and so there exists a
standard $A_n$-module $\Delta_n(\lambda)$ for each weight $\lambda$ in
$\Lambda_n$, such that the associated simple $L_n(\lambda)$ arises as
its head. The restriction of a standard module has a filtration by
standards and we denote by $\supp_n(\lambda)$ the multi-set of labels
for standard modules occurring in such a filtration of
$\res_n\Delta_n(\lambda)$. The embedding of $\Lambda_n$ in
$\Lambda_{n+2}$ induces an embedding of $\supp_n(\lambda)$ inside
$\supp_{n+2}(\lambda)$, which becomes an identification if
$\lambda\in\Lambda_{n-2}$. We denote by $\supp(\lambda)$
the set $\supp_n(\lambda)$ with $n>>0$.

Suppose that we have determined the blocks of such a family of
algebras (or at least a linkage principle); we are thinking of the
cases where we have an alcove geometry at hand, but will avoid stating
the result in that form. Let $\res_n^{\lambda}$ be the functor
$\pr_{n-1}^{\lambda}\res_n$ and $\ind_n^{\lambda}$ be the functor
$\pr_{n+1}^{\lambda}\ind_n$ where $\pr_{n}^{\lambda}$ is projection
onto the block containing $\lambda$ for $A_{n}$.  Note that
$\res_n^{\lambda}$ is exact and $\ind_n^{\lambda}$ is right exact.
We will regard these functors as analogues of \emph{translation
  functors} in Lie theory.

Let ${\mathcal B}_n(\lambda)$ denote the set of weights in the block
of $A_n$ which contains $\lambda$.  Our embedding of $\Lambda_n$ into
$\Lambda_{n+2}$ induces an embedding of ${\mathcal B}_n(\lambda)$ into
${\mathcal B}_{n+2}(\lambda)$, and we denote by ${\mathcal
  B}(\lambda)$ the corresponding limiting set. We will say that two
elements $\lambda$ and $\lambda'$ are \emph{translation equivalent} 
if they satisfy 
\begin{enumerate}
\item[(i)] The weight $\lambda'$ is the only element of ${\mathcal
  B}(\lambda')\cap\supp(\lambda)$.
\item[(ii)] The weight $\lambda$ is the only element of ${\mathcal
  B}(\lambda)\cap\supp(\lambda')$.
\item[(iii)] For all weights $\mu\in {\mathcal B}(\lambda)$ there is a unique
  element $\mu'\in {\mathcal B}(\lambda')\cap \supp(\mu)$, and
  $\mu$ is the unique element in ${\mathcal B}(\lambda)\cap
  \supp(\mu')$.
\end{enumerate}
Clearly conditions (i) and (ii) are special cases of (iii); we list
them separately as in an alcove geometry (i) and (ii) will be enough
for (iii) to hold. When $\lambda$ and $\lambda'$ are translation
equivalent then we will denote by $\theta:{\mathcal
  B}(\lambda)\rightarrow{\mathcal B}(\lambda')$ the bijection taking
$\mu$ to $\mu'$. We will see that translation equivalent weights
belong to Morita equivalent blocks.

We will put a very crude partial order on weights in ${\mathcal
  B}(\lambda)$ by saying that $\lambda>\mu$ if there exists $n$ such
that $\mu\in\Lambda_n$ and $\lambda\in\Lambda_{n+2t}$ for some $t\in
\NN$, but $\lambda\notin\Lambda_n$. Note that this is the opposite of
the standard order arising from the quasi-hereditary structure; we
prefer to work with the natural order on the size of partitions.  In
the following proposition, by a unique element in a multi-set we mean
one with multiplicity one.

\begin{prop}\label{transmain}
Let $A$ be a tower of recollement. Suppose that $\lambda\in\Lambda_n$
and $\lambda'\in\Lambda_{n-1}$ are translation equivalent, and that
$\mu\in {\mathcal B}_n(\lambda)$ is such that the $\mu'$ is in
${\mathcal B}_{n-1}(\lambda')$.  Then we have that
\begin{equation}\label{simres}
\res_n^{\lambda'}L_n(\mu)\cong L_{n-1}(\mu')\quad\quad\mbox{
  and}\quad\quad
 \ind_{n-1}^{\lambda}L_{n-1}(\mu')\cong L_{n}(\mu)
\end{equation}
 for all $\mu\in{\mathcal B}_n(\lambda)$. Further if $\tau\in
 {\mathcal B}_n(\lambda)$ is such that $\tau'$ is in
 ${\mathcal B}_{n-1}(\lambda')$ then we have 
\begin{equation}\label{decres}
[\Delta_n(\mu):L_n(\tau)]=[\Delta_{n-1}(\mu'):L_{n-1}(\tau')]
\end{equation}
and
\begin{equation}\label{homsame}
\Hom(\Delta_n(\mu),\Delta_n(\tau))\cong
\Hom(\Delta_{n-1}(\mu'),\Delta_{n-1}(\tau')).
\end{equation} 
\end{prop}
\begin{proof}
We begin with (\ref{simres}). Consider the exact sequence
$$\Delta_n(\mu)\too L_n(\mu)\too 0.$$ Applying $\res_n^{\lambda'}$ we
obtain by our assumptions the exact sequence
$$\Delta_{n-1}(\mu')\too \res_n^{\lambda'}L_n(\mu)\too 0$$
and hence $\res_n^{\lambda'}L_n(\mu)$ has simple head
$L_{n-1}(\mu')$, and possibly other composition factors
$L_{n-1}(\tau')$ with $\tau'>\mu'$.

If $L(\tau')$ is in the socle of $\res_n^{\lambda'}L_n(\mu)$ then
we have 
\begin{eqnarray*}
\Hom(\Delta_n(\tau),L_n(\mu))=&&\!\!\!\!\!\!
\Hom(\ind_{n-1}^{\lambda}\Delta_{n-1}(\tau'),L_n(\mu))\cong
\Hom(\ind_{n-1}\Delta_{n-1}(\tau'),L_n(\mu))\\
&\cong&
\Hom(\Delta_{n-1}(\tau'),\res_nL_n(\mu))\cong
\Hom(\Delta_{n-1}(\tau'),\res_n^{\lambda'}L_n(\mu))
\neq 0
\end{eqnarray*} by our
assumptions and Frobenius reciprocity, which implies that $\tau=\mu$
and hence $\tau'=\mu'$. Therefore $\res_n^{\lambda'}L_n(\mu)\cong
L_{n-1}(\mu')$  as required.

Next consider the exact sequence
$$\Delta_{n-1}(\mu')\too L_{n-1}(\mu')\too 0.$$ Applying
$\ind_{n-1}^{\lambda}$ we obtain by our assumptions the exact sequence
$$\Delta_{n}(\mu)\too \ind_{n-1}^{\lambda}L_{n-1}(\mu')\too 0.$$
and hence $\ind_{n-1}^{\lambda}L_{n-1}(\mu')$ has simple head
$L_{n}(\mu)$, and possibly other composition factors
$L_{n}(\tau)$ with $\tau>\mu$. Now apply $\res_n^{\lambda'}$ to
obtain the exact sequence
$$\Delta_{n-1}(\mu')\too
\res_n^{\lambda'}\ind_{n-1}^{\lambda}L_{n-1}(\mu')\too 0.$$ Then
$\res_n^{\lambda'}\ind_{n-1}^{\lambda}L_{n-1}(\mu')$ has simple head
$L_{n-1}(\mu')$ and possibly other composition factors
$L_{n-1}(\tau')$ with $\tau'>\mu'$ corresponding to those in
$\ind_{n-1}^{\lambda}L_{n-1}(\mu')$. We have
$$\Hom(L_{n-1}(\mu'),\res_n^{\lambda'}\ind_{n-1}^{\lambda}L_{n-1}(\mu'))
\cong
\Hom(\ind_{n-1}^{\lambda}L_{n-1}(\mu'),\ind_{n-1}^{\lambda}L_{n-1}(\mu'))
\neq 0$$
and hence $L_{n-1}(\mu')$ must appear in the socle of 
$\res_n^{\lambda'}\ind_{n-1}^{\lambda}L_{n-1}(\mu')$. This forces
$$\res_n^{\lambda'}\ind_{n-1}^{\lambda}L_{n-1}(\mu')\cong
L_{n-1}(\mu')$$
and as we already have that $\res_n^{\lambda'}$ takes simples to
simples we deduce that
$$\ind_{n-1}^{\lambda}L_{n-1}(\mu')\cong
L_{n}(\mu)$$
which completes our proof of (\ref{simres}). Now (\ref{decres})
follows immediately as $\res_{n}^{\lambda'}\Delta_n(\mu)\cong
\Delta_{n-1}(\mu')$, while (\ref{homsame}) follows from
\begin{eqnarray*}
\Hom(\Delta_n(\mu),\Delta_n(\tau))
&\cong&\Hom(\ind^{\lambda}_{n-1}\Delta_{n-1}(\mu'),\Delta_n(\tau))\\
&\cong &
\Hom(\Delta_{n-1}(\mu'),\res_n^{\lambda'}\Delta_n(\tau))
\cong \Hom(\Delta_{n-1}(\mu'),\Delta_{n-1}(\tau')).
\end{eqnarray*}
\end{proof}

Let $P_n(\lambda)$ denote the projective cover of $L_n(\lambda)$. As
our algebras are quasihereditary we have that $P_n(\lambda)$ has a
filtration by standard modules with well-defined filtration
multiplicities; we denote the multiplicity of $\Delta_n(\mu)$ in such a
filtration by $(P_n(\lambda):\Delta_n(\mu))$.

\begin{prop}\label{transproj}
Suppose that $\lambda\in\Lambda_n$ and $\lambda'\in\Lambda_{n-1}$ are
translation equivalent. Then for all $\mu\in{\mathcal B}_n(\lambda)$
with $\mu'\in{\mathcal B}_{n-1}(\lambda')$ we have
\begin{equation}\label{indproj}
\ind_{n-1}^{\lambda}P_{n-1}(\mu')\cong P_n(\mu).
\end{equation}
If $\mu\in{\mathcal B}_{n-2}(\lambda)$ we have
\begin{equation}\label{resproj}
\res_{n}^{\lambda'}P_{n}(\mu)\cong P_{n-1}(\mu').
\end{equation}
\end{prop}
\begin{proof} We begin with (\ref{indproj}).
The functor $\ind_{n-1}$ takes projectives to projectives, and hence
so does $\ind_{n-1}^{\lambda}$. We must show that inducing an
indecomposable projective gives an indecomposable projective with the
right weight.

Suppose we have an exact sequence
$$\ind_{n-1}^{\lambda}P_{n-1}(\mu')\rightarrow L_n(\tau)\rightarrow
0$$
for some $\tau\in{\mathcal B}_n(\lambda)$. Then we have
\begin{eqnarray*}
0&\neq &
\Hom_n(\ind_{n-1}^{\lambda}P_{n-1}(\mu'),L_n(\tau))\\
&\cong &
\Hom_{n-1}(P_{n-1}(\mu'),\res_n^{\lambda'}L_n(\tau))
\cong \Hom_{n-1}(P_{n-1}(\mu'),L_n(\tau'))
\end{eqnarray*}
by Proposition \ref{transmain}. Therefore we must have $\mu'=\tau'$,
and hence $\mu=\tau$ and
$$\Hom_n(\ind_{n-1}^{\lambda}P_{n-1}(\mu'),L_n(\mu))\cong k.$$
This implies that $\ind_{n-1}^{\lambda}P_{n-1}(\mu')$ has simple head
$L_n(\mu)$ and hence is isomorphic to $P_n(\mu)$.

Next we consider (\ref{resproj}). As $A_ne_n$ is a direct summand of
the left $A_n$-module $A_n$, it is a projective
$A_n$-module. Moreover, as $e_nA_ne_n\cong A_{n-2}$ we have that $A_ne_n$
contains precisely those indecomposable projective $A_n$-modules
labelled by weights in $\Lambda_{n-2}$. In a tower of recollement we
have that 
$$\res_{n-1}A_ne_n\cong A_{n-1}$$ 
as a left $A_{n-1}$-module. This implies that for
$\mu\in\Lambda_{n-2}$, the module $\res_{n}P_n(\mu)$ (and hence 
$\res_{n}^{\lambda'}P_n(\mu)$)
is projective.

As $\res_n^{\lambda'}$ is an exact functor and $P_n(\mu)$ has
simple head $L_n(\mu)$ we know from Proposition \ref{transmain} 
that
\begin{equation}\label{fakesum}
\res_n^{\lambda'}P_n(\mu)=P_n(\mu')\oplus Q
\end{equation}
for some projective $A_{n-1}$-module $Q$.  However, by
Brauer-Humphreys reciprocity and Proposition \ref{transmain} we have
$$(P_n(\mu):\Delta_n(\tau))=[\Delta_n(\tau):L_n(\mu)]=
[\Delta_{n-1}(\tau'):L_{n-1}(\mu')]
=(P_{n-1}(\mu'):\Delta_{n-1}(\tau')).$$ As $\res_n^{\lambda'}$ is
exact and takes $\Delta_n(\tau)$ to $\Delta_{n-1}(\tau')$ this implies
that $Q=0$.
\end{proof}

We would like to argue that two blocks labelled by translation
equivalent weights are Morita equivalent. However, the fact that not
every projective module restricts to a projective in (\ref{resproj})
causes certain complications. 

\begin{lem} \label{Pstable}
If $\lambda\in\Lambda_n$ then
$$G_n(P_n(\lambda))\cong P_{n+2}(\lambda).$$
\end{lem}
\begin{proof}
By \cite[Chapter I, Theorem 6.8]{ass1} $G_n(P_n(\lambda))$ is an
indecomposable projective. We have an exact sequence
$$P_n(\lambda)\rightarrow \Delta_n(\lambda)\rightarrow 0$$ And hence
as $G_n$ is right exact and takes standards to standards we obtain
$$G_n(P_n(\lambda))\rightarrow \Delta_{n+2}(\lambda)\rightarrow 0.$$
This implies that $G_n(P_n(\lambda))\cong P_{n+2}(\lambda)$.
\end{proof}

\begin{lem}\label{Phomstable}
If $\mu,\tau\in\Lambda_n$ then
$$\Hom_n(P_n(\mu),P_n(\tau))\cong\Hom_{n+2}(P_{n+2}(\mu),P_{n+2}(\tau))$$
and this extends to an algebra isomorphism
$$\End_n(\bigoplus_{\mu\in\Gamma}P_n(\mu))\cong
\End_{n+2}(\bigoplus_{\mu\in\Gamma}P_{n+2}(\mu))$$
where $\Gamma\subseteq \Lambda_n$.
\end{lem}
\begin{proof}
See \cite[Chapter I, Theorem 6.8]{ass1}.
\end{proof}

The algebra $A_n$ decomposes as a direct sum of
indecomposable projective modules:
$$A_n=\bigoplus_{\lambda\in\Lambda_n}P_n(\lambda)^{d_{n,\lambda}}$$
for some integers $d_{n,\lambda}$. There is a corresponding
decomposition of $1\in A_n$ as a sum of (not necessarily primitive)
orthogonal idempotents
$1=\sum_{\lambda\in\Lambda_n}e_{n,\lambda}$ where
$A_ne_{n,\lambda}=P_n(\lambda)^{d_{n,\lambda}}$. As $\Lambda_n$
decomposes as a union of blocks the algebra $A_n$ decomposes as a
direct sum of (block) subalgebras
$$A_n=\bigoplus_{\lambda}A_n(\lambda)$$
where the sum runs over a set of block representatives and 
$$A_n(\lambda)=\bigoplus_{\mu\in{\mathcal B}_n(\lambda)}P_n(\mu)^{d_{n,\mu}}.$$

Now let $\Gamma\subset {\mathcal B}_n(\lambda)$ and consider the idempotent
$e_{n,\Gamma}=\sum_{\gamma\in\Gamma}e_{n,\gamma}$. We define the algebra
$A_{n,\Gamma}(\lambda)$ by
$$A_{n,\Gamma}(\lambda)=e_{n,\Gamma}A_n(\lambda)e_{n,\Gamma}.$$ 
By Lemma \ref{Phomstable} we have that $A_{n,\Gamma}(\lambda)$ and
$A_{m,\Gamma}(\lambda)$ are Morita equivalent for all $m$ such that
$\Gamma\subset {\mathcal B}_m(\lambda)$. 

\begin{thm}\label{genmor}
Suppose that $\lambda$ and $\lambda'$ are translation equivalent with
$\lambda\in \Lambda_n$, and set 
$$\Gamma=\theta({\mathcal B}_n(\lambda))\subset {\mathcal
  B}_{n+1}(\lambda').$$ Then $A_n(\lambda)$ and
$A_{n+1,\Gamma}(\lambda')$ are Morita equivalent.  In particular, if
there exists an $n$ such that $|{\mathcal B}_n(\lambda)|=|{\mathcal
  B}_{n+1}(\lambda')|$ then $A_n(\lambda)$ and $A_{n+1}(\lambda')$ are
Morita equivalent.
\end{thm}
\begin{proof}

We will show that the basic algebras corresponding to $A_n(\lambda)$
and $A_{n+1,\Gamma}(\lambda')$ are isomorphic; i.e. that
$$\End_n(\bigoplus_{\mu\in{\mathcal B}_n(\lambda)}P_n(\mu))\cong
\End_{n+1}(\bigoplus_{\nu'\in\Gamma}P_{n+1}(\nu')).$$
By Lemma \ref{Phomstable} it is enough to show that 
$$\End_{n+2}(\bigoplus_{\mu\in{\mathcal B}_n(\lambda)}P_{n+2}(\mu))\cong
\End_{n+1}(\bigoplus_{\nu'\in\Gamma}P_{n+1}(\nu')).$$

Suppose that $\mu,\tau\in{\mathcal B}_n(\lambda)$. Then by Lemma
\ref{Phomstable} and Proposition \ref{transproj} we have
\begin{eqnarray*}
\Hom_{n+2}(P_{n+2}(\mu),P_{n+2}(\tau))
&\cong&\Hom_{n+2}(\ind_{n+1}^{\lambda}P_{n+1}(\mu'),P_{n+2}(\tau))\\ 
&\cong&\Hom_{n+1}(P_{n+1}(\mu'),\res_{n+2}^{\lambda'}P_{n+2}(\tau))\\ 
&\cong&\Hom_{n+1}(P_{n+1}(\mu'),P_{n+1}(\tau'))
\end{eqnarray*} 
Next we will show that these isomorphisms are also
compatible with the multiplicative structure in each of our algebras.

Let $P$, $Q$, and $R$ be indecomposable projectives for $A_{n+2}$
labelled by elements from ${\mathcal B}_n(\lambda)$. Then there exist
indecomposable projectives $P'$, $Q'$ and $R'$ labelled by elements in
${\mathcal B}_{n+1}(\lambda')$ such that
$$P=\ind_{n+1}^{\lambda}P'\quad\quad Q=\ind_{n+1}^{\lambda}Q'\quad\quad
Q'=\res_{n+2}^{\lambda'}Q\quad\quad R'=\res_{n+2}^{\lambda'}R.$$
An isomorphism $\alpha$ giving a Frobenius reciprocity of the form
$$\Hom_{n+1}(M,\res_{n+2}N)\cong \Hom_{n+2}(\ind_{n+1}M,N)$$
is given by the map taking $\phi$ to $\alpha(\phi)$ where 
$$\alpha(\phi)(a\otimes m)=a\phi(m)$$ for all $a\in A_{n+2}$ and $m\in M$,
and extending by linearity. (Recall that $\ind_{n+1}$ is just the
function $A_{n+2}\otimes_{A_{n+1}}\!\!$--.) Given 

\begin{equation}\label{comdiag}
\xymatrix@R=15pt@C=10pt{\Hom_{n+2}(P,Q) \ar@{}[r]|{\times}\ar[d]_{\cong}
&\Hom_{n+2}(Q,R)\ar[r]\ar[d]_{\cong}&\Hom_{n+2}(P,R)\ar[d]_{\cong}\\
\!\!\!\!\!\!\!\!\!\Hom_{n+2}(\ind_{n+1}^{\lambda}P',\ind_{n+1}^{\lambda}Q')\ \ 
\ar@{}[r]|{\times}&
\!\!\!\!\!\!\Hom_{n+2}(\ind_{n+1}^{\lambda}Q',R)\ar[r]&
\Hom_{n+2}(\ind_{n+1}^{\lambda}P',R)\\
\!\!\!\!\!\!\!\!\!
\Hom_{n+1}(P',\res_{n+2}^{\lambda'}\ind_{n+1}^{\lambda}Q'))\ \ 
\ar@{}[r]|{\times}
\ar[u]^{\alpha}&
\!\!\!\!\!\!\Hom_{n+1}(Q'\res_{n+2}^{\lambda'}R)\ar[r]\ar[u]^{\alpha}&
\Hom_{n+1}(P'\res_{n+2}^{\lambda'}R)\ar[u]^{\alpha}\\
\Hom_{n+1}(P',Q')
\ar@{}[r]|{\times}\ar@{}[d]|{\phi}\ar[u]^{\cong}
&\Hom_{n+1}(Q',R')\ar@{}[d]|{\psi}\ar[u]^{\cong}\ar[r]
&\Hom_{n+1}(P',R')\ar@{}[d]|{\psi\circ\phi}\ar[u]^{\cong}\\
&&}
\end{equation}
we need to check that
$\alpha(\psi\circ\phi)=\alpha(\psi)\circ\alpha(\phi)$.
We have 
$$\alpha(\phi)(\sum_ia_i\otimes p_i)=\sum_ia_i\phi(p_i)$$
where $a_i\in A_{n+2}$ and $p_i\in P'$. As $\phi(p_i)\in Q'
\cong\res_{n+2}^{\lambda'}(\ind_{n+1}^{\lambda}Q')$ we have
$$\phi(p_i)=\sum_ja'_j\otimes q_j$$
where $a'_j\in A_{n+2}$ and $q_j\in Q'$. Now 
\begin{eqnarray*}
(\alpha(\psi)\circ\alpha(\phi))(\sum_i a_i\otimes p_i)&=&
\alpha(\psi)\Bigl(\sum_ia_i\bigl(\sum_ja'_j\otimes q_j\bigr)\Bigr)\\
&=&\alpha(\psi)\Bigl(\sum_{i,j}a_ia'_j\otimes q_j\Bigr)
=\sum_{i,j}a_ia'_j\psi(q_j)
\end{eqnarray*}
where the second equality follows from the action of $A_{n+2}$ on
$\ind_{n-1}Q'$. On the other hand
\begin{eqnarray*}
\alpha(\psi\circ\phi)(\sum_i a_i\otimes p_i)&=&
\sum_ia_i(\psi\circ\phi)(p_i)\\
&=&\sum_ia_i\psi(\sum_ja'_j\otimes q_j)\\
&=&\sum_ia_i\sum_ja_j'\psi(q_j)
=\sum_{i,j}a_ia_j'\psi(q_j)
\end{eqnarray*}
and hence $\alpha(\psi\circ\phi)=\alpha(\psi)\circ\alpha(\phi)$ as required.
\end{proof}

Next we will consider how we can relate the cohomology of
$A_{n+1,\Gamma}(\lambda')$ to that of $A_{n+1}(\lambda')$ and hence
compare the cohomology of $A_n(\lambda)$ with that of
$A_{n+1}(\lambda')$. We say that a subset $\Gamma\subset {\mathcal
  B}_n(\lambda)$ 
is \emph{saturated} if $\mu\in\Gamma$ and
$\nu\in{\mathcal B}_n(\lambda)$ with $\nu>\mu$ implies that
$\nu\in\Gamma$. A subset $\Gamma\subset {\mathcal B}_n(\lambda)$ is
\emph{cosaturated} if ${\mathcal B}_n(\lambda)\backslash\Gamma$ is
saturated.

\begin{lem}\label{cosat}
The set $\Gamma=\theta({\mathcal B}_n(\lambda))$ is cosaturated in
${\mathcal B}_{n+1}(\lambda')$.
\end{lem}
\begin{proof}
We need to show that if $\mu'\in{\mathcal B}_{n+1}(\lambda')\backslash
\Gamma$ and $\nu'\in{\mathcal B}_{n+1}(\lambda')$ with $\nu'>\mu'$
then $\nu'\in{\mathcal B}_{n+1}(\lambda')\backslash \Gamma$. Suppose
for a contradiction that $\nu'\in\Gamma$. Then
$\nu=\theta^{-1}(\nu')\in{\mathcal B}_n(\lambda)$ and
$\mu=\theta^{-1}(\mu')\notin{\mathcal B}_n(\lambda)$. As
$\mu\in{\mathcal B}(\lambda)$ we must have $|\mu|\geq n+2$, and as
$\mu'\in\supp(\mu)$ we have 
$$|\mu'|\geq n+2\pm 1\geq n+1.$$
Now $|\nu|\leq n$ and so 
$$|\nu'|\leq n\pm 1\leq n+1$$
but this contradicts the assumption that $\mu'<\nu'$.
\end{proof}

If $\Gamma\subset {\mathcal B}_n(\lambda)$ is cosaturated then by
\cite[A.3.11]{don2} the algebra $A_{n,\Gamma}(\lambda)$ is
quasihereditary with standard modules given by
$$\{e_{n,\Gamma}\Delta_n(\mu):\mu\in\Gamma\}.$$
Moreover if $X$ is any $A_n$-module having a $\Delta$-filtration with
factors $\Delta_n(\mu)$ for $\mu\in\Gamma$, and $Y$ is any $A_n$-module,
then for all $i\geq 0$ we have \cite[A.3.13]{don2}
$$\Ext^i_n(X,Y)=\Ext^i_{A_n(\lambda)}(X,Y)\cong
\Ext^i_{A_{n,\Gamma}(\lambda)}(e_{n,\Gamma}X, e_{n,\Gamma}Y).$$

Combining the above remarks with Theorem \ref{genmor} and Lemma \ref{cosat}
we obtain

\begin{cor}\label{extsame}
If $\lambda\in\Lambda_n$ and $\lambda'\in\Lambda_{n+1}$ are
translation equivalent then for all $i\geq 0$ and for all
$\mu\in{\mathcal B}_n(\lambda)$ we have
$$\Ext^i_n(\Delta_n(\lambda),\Delta_n(\mu))\cong
\Ext^i_{n+1}(\Delta_{n+1}(\lambda'),\Delta_{n+1}(\mu')).$$
\end{cor}

We will say that two weights $\lambda$ and $\lambda'$ are in the same
\emph{translation class} if they are related by the equivalence
relation generated by translation equivalence.  Then analogues of
(\ref{decres}), (\ref{homsame}), Theorem \ref{genmor} and Corollary
\ref{extsame} also hold for weights in the same translation class.

In Lie theory one can consider translation between two weights in the
same facet (corresponding to the case considered above), or from one
facet to another. We will now give an analogue of translation onto a wall
for a tower of recollement.

We will say that $\lambda'$ \emph{separates} $\lambda^-$ and
$\lambda^+$ if 
\begin{enumerate}
\item[(i)] The weight $\lambda'$ is the only element of ${\mathcal
  B}(\lambda')\cap\supp(\lambda^-)$.
\item[(ii)] The weight $\lambda'$ is the only element of ${\mathcal
  B}(\lambda')\cap\supp(\lambda^+)$.
\item[(iii)] The weights $\lambda^+$ and $\lambda^-$ are the only
  elements of ${\mathcal
  B}(\lambda^-)\cap\supp(\lambda')$.
\end{enumerate}
Whenever we consider a pair of weights $\lambda^-$ and $\lambda^+$
separated by $\lambda'$ we shall always assume that
$\lambda^-<\lambda^+$.

\begin{thm}\label{ontowall}
(i) If $\lambda'\in\Lambda_{n-1}$ separates $\lambda^-$ and
  $\lambda^+$ then
$$\res^{\lambda'}_nL_n(\lambda^+)\cong L_{n-1}(\lambda').$$
(ii) If further we have
  $\Hom(\Delta_n(\lambda^+),\Delta_n(\lambda^-))\neq 0$ then
$$\res^{\lambda'}_nL_n(\lambda^-)=0$$
and $\ind^{\lambda^-}_{n-1}\Delta_n(\lambda')$ is a nonsplit extension of
$\Delta_n(\lambda^-)$ by $\Delta_n(\lambda^+)$ and has simple head
  $L_n(\lambda^+)$. 
\end{thm}
\begin{proof}
Arguing as in the proof of Proposition \ref{transmain} we see that
$\res_n^{\lambda'}L_n(\lambda^i)$ is either $0$ or has simple head
$L_{n-1}(\lambda')$ for $i=+,-$. Also, any other composition factors
$L_{n-1}(\tau')$ of $\res_n^{\lambda'}L_n(\lambda^i)$ must satisfy
$\tau'>\lambda'$. Note that by assumption we have a short exact
sequence
\begin{equation}\label{sepses}
0\too
\Delta_n(\lambda^-)\too\ind_{n-1}^{\lambda^-}\Delta_{n-1}(\lambda')
\too \Delta_n(\lambda^+)\too 0
\end{equation}
and hence
\begin{equation}\label{twohom}
\Hom(\Delta_{n-1}(\tau'),\res_n^{\lambda'}L_n(\lambda^i))\cong
\Hom(\ind_{n-1}^{\lambda^-}\Delta_{n-1}(\tau'),L_n(\lambda^i))
\end{equation}
is non-zero when $i=+$ and $\tau'=\lambda'$, and is zero when
$\tau'>\lambda'$ by our assumptions. This completes the proof of (i).

Now suppose that $\Hom(\Delta_n(\lambda^+),\Delta_n(\lambda^-))\neq
0$. Then we have that
$$[\Delta_n(\lambda^-):L_n(\lambda^+)]\neq 0.$$
By exactness and the first part of the Theorem, the unique copy of
$L_{n-1}(\lambda')$ in
$$\res_{n}^{\lambda'}\Delta_n(\lambda^-)\cong\Delta_{n-1}(\lambda')$$
must come from $\res_n^{\lambda'}L_n(\lambda^+)$, and hence
$\res_n^{\lambda'}L_n(\lambda^-)$ cannot have simple head
$L_{n-1}(\lambda')$. But this implies by the first part of the proof
that $\res_n^{\lambda'}L_n(\lambda^-)=0$. Therefore the Hom-space in
(\ref{twohom}) must be zero when $\tau'=\lambda'$ and $i=1$, which
implies that (\ref{sepses}) is a non-split extension whose central
module has simple head $L_n(\lambda^+)$ as required.
\end{proof}

Suppose that $\lambda'$ and $\lambda^+$ are weights with
$\lambda'<\lambda^+$ and $\lambda'\in\supp(\lambda^+)$ such that for
every weight $\tau'\in{\mathcal B}(\lambda')$ either (i) there is a unique
weight $\tau^+\in{\mathcal B}(\lambda^+)\cap\supp(\tau')$ and $\tau'$
is the unique weight in ${\mathcal B}(\lambda')\cap\supp(\tau^+)$, or
(ii) there exists $\tau^-,\tau^+\in{\mathcal B}(\lambda^+)$ such that
$\tau'$ separates $\tau^-$ and $\tau^+$. Then we say that $\lambda'$
is in the \emph{lower closure} of $\lambda^+$. If further 
$$\Hom(\Delta_n(\tau^+),\Delta_n(\tau^-))\neq 0$$ whenever $\tau'
\in{\mathcal B}(\lambda')$
separates $\tau^-$ and $\tau^+$ in ${\mathcal B}(\lambda^+)$ then we
shall say that ${\mathcal B}(\lambda^+)$ has \emph{local
  homomorphisms} with respect to ${\mathcal B}(\lambda')$.

\begin{prop}\label{alcok}
Suppose that $\lambda'\in\Lambda_{n-1}$ is in the lower closure of
$\lambda^+\in\Lambda_n$, and that ${\mathcal B}(\lambda^+)$ has enough
local homomorphisms with respect to ${\mathcal B}(\lambda')$. Then
$$[\Delta_{n-1}(\lambda'):L_{n-1}(\mu')]=[\Delta_n(\lambda^+):L_n(\mu^+)].$$
\end{prop}
\begin{proof}
We have by our assumptions that 
$$\res_n^{\lambda'}\Delta_n(\lambda^+)\cong \Delta_{n-1}(\lambda').$$
As $\res_n^{\lambda'}$ is an exact functor, it is enough to determine
its effect on simples $L_n(\mu^+)$ in $\Delta_n(\lambda^+)$. If there
exists $\mu'$ separating $\mu^+$ from $\mu^-$ then the result follows
from Theorem \ref{ontowall}, while if $\mu^+$ is the only element
in ${\mathcal B}(\lambda^+)\cap\supp(\mu')$ then if follows as in the
proof of Proposition \ref{transmain}.
\end{proof}

Thus, as long as there are enough local homomorphisms, the
decomposition numbers for $\Delta_n(\lambda)$ determine those for all
weights in the lower closure of $\lambda$.

We can generalise the results of this section up to Corollary \ref{extsame} 
by replacing
$\res^{\lambda}_n$ and $\ind^{\lambda}_n$ by any pair of functor
families $R_n$ and $I_n$ with the following properties.

Let $A$ be a tower of recollement, with
$\lambda,\lambda'\in\Lambda$, and fix $i\in\NN$. Suppose that
we have functors
$$R_n:A_n\mbox{\rm -mod}\rightarrow A_{n-i}\mbox{\rm -mod}$$ 
for $n\geq i$ and
$$I_n:A_n\mbox{\rm -mod}\rightarrow A_{n+i}\mbox{\rm -mod}$$ 
for $n\geq 0$ satisfying
\begin{enumerate}
\item[(i)]
The functor $I_n$ is left adjoint to $R_{n+i}$ for all $n$.
\item[(ii)]
The functor $R_n$ is exact and $I_n$ is right exact for all $n$ where
they are defined.
\item[(iii)]
There is a bijection $\theta:{\mathcal B}(\lambda)\rightarrow {\mathcal
  B}(\lambda')$ taking $\mu$ to $\mu'$ such that for all
$n\geq i$, if $\mu\in{\mathcal B}_n(\lambda)$ and $\mu'\in{\mathcal
  B}_{n-i}(\lambda')$ then 
$$R_n\Delta_n(\mu)\cong \Delta_{n-i}(\mu')\quad\quad\wand\quad\quad
I_{n-i}\Delta_{n-i}(\mu')=\Delta_n(\mu)$$
and $R_n\Delta_n(\mu)=0$ otherwise.
\item[(iv)] If $\Gamma_{n}=\theta({\mathcal
  B}_n(\lambda))\subset{\mathcal B}_m(\lambda')$ for some $m$ then
  $\Gamma_n$ is cosaturated in ${\mathcal B}_{m}(\lambda')$.
\item[(v)] There exists $t\in\NN$ such that for all $n$ and for all
  $\mu\in{\mathcal B}_{n-t}(\lambda)$ the module $R_nP_n(\mu)$ is
  projective. 
\end{enumerate}
Then we say that $\lambda$ and $\lambda'$ are
\emph{$(R,I)$-translation equivalent}. In this case the proofs of
Propositions \ref{transmain} and \ref{transproj} go through
essentially unchanged, and we get 

\begin{thm} \label{ripair}
Suppose that $\lambda$ and $\lambda'$ are
  $(R,I)$-translation equivalent and $n\geq i$. Then for all $\mu\in{\mathcal
    B}_n(\lambda)$ with $\mu'\in{\mathcal B}_{n-i}(\lambda')$
 we have
$$R_nL_n(\mu)\cong L_{n-i}(\mu'),\quad\quad\quad\quad
I_{n-i}L_{n-i}(\mu')\cong L_n(\mu)\quad\quad\wand\quad\quad
I_{n-i}P_{n-i}(\mu')\cong P_n(\mu)$$
and if $\mu\in{\mathcal B}_{n-t}(\lambda)$ then
$$R_nP_n(\mu)\cong P_{n-i}(\mu').$$
Moreover, if the adjointness isomorphism 
$$\alpha: \Hom_{n-i}(M,R_n(N))\rightarrow \Hom_n(I_{n-i}(M),N)$$ is
multiplicative (i.e. makes the diagram (\ref{comdiag}) commute) then
there is a Morita equivalence between $A_n(\lambda)$ and
$A_{n+i,\Gamma_{n+i}}(\lambda')$ and for all $\mu,\tau\in{\mathcal
  B}_n(\lambda)$ and $j\geq 0$ we have
$$\Ext^j_{n}(\Delta_n(\mu),\Delta_n(\tau))\cong
\Ext^j_{n+i}(\Delta_{n+i}(\mu'),\Delta_{n+i}(\tau')).$$
\end{thm}

\section{A generalised restriction/induction pair}\label{genir}

We wish to show (in Section \ref{brauertran}) that two weights in the
same facet for the Brauer algebra give rise to Morita equivalent
blocks (at least when we truncate the blocks to have the same number
of simples). However, the usual induction and restriction functors are
not sufficient to show this except in the alcove case. To remedy this,
in this section we will consider a variation on the usual induction
and restriction functors. As $\delta$ will be fixed throughout, we
will denote $B_n(\delta)$ simply as $B_n$.

\begin{figure}[ht]
\center{\includegraphics{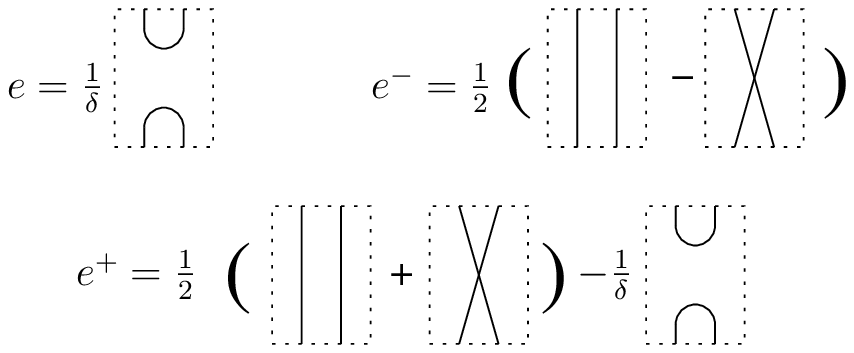}}
\caption{Idempotents in $B_2$}
\label{idemsare}
\end{figure}

First consider $B_2$ with $\delta\neq 0$. It is easy to see
that this is a semisimple algebra, with a decomposition
$$1=e+e^-+e^+$$
of the identity into primitive orthogonal idempotents given by
the elements in Figure \ref{idemsare}.

There are three standard modules for this algebra, which we will
denote by
\begin{eqnarray*}
&\Delta_2(0)&=\,<e>\\
S^-=&\Delta_2(1,1)&=\,<e^->\\
S^+=&\Delta_2(2)&=\,<e^+>
\end{eqnarray*}

For $n\geq 2$ consider the subalgebra $B_{n-2}\otimes
B_2\subseteq B_n$ obtained by letting
$B_{n-2}$ act on the leftmost $n-2$ lines and $B_2$
act on the rightmost pair of lines. We will view elements of
$B_{n-2}$ and $B_2$ as elements of $B_n$ via
this embedding. Note that under this embedding the two algebras
obviously commute with each other. 

In particular, for any $B_n$-module $M$ the vector spaces
$e^{\pm}M$ are $B_{n-2}$-modules. Thus we have a pair of functors
$\res_n^{\pm}$ from $B_n$-mod to $B_{n-2}$-mod given on
objects by the map $M\longmapsto e^{\pm}M$. Note that these functors can
also be defined as 
$$\res_n^{\pm}M=e^{\pm}\res_{B_{n-2}\otimes B_2}^{B_n}M.$$

We have
\begin{eqnarray*}
&\Hom_{B_n}(\ind_{B_{n-2}\otimes B_2}^{B_n}(N\boxtimes S^{\pm}),M)&
\cong \Hom_{B_{n-2}\otimes B_2}(N\boxtimes S^{\pm},
\res_{B_{n-2}\otimes B_2}^{B_n}M)\\
&&\cong \Hom_{B_{n-2}}(N,e^{\pm}\res_{B_{n-2}\otimes  B_2}^{B_n}M)
\end{eqnarray*}
and so the functors $\ind_{n-2}^{\pm}$ from $B_{n-2}$-mod to
$B_n$-mod given by 
$$\ind_{n-2}^{\pm}N=\ind_{B_{n-2}\otimes B_2}^{B_n}(N\boxtimes S^{\pm})$$
are left adjoint to $\res_n^{\pm}$.

\begin{lem} \label{altdef}
Let $N$ be a $B_{n-2}$-module. Then we have
$$\ind_{n-2}^{\pm}N\cong B_ne^{\pm}\otimes_{B_{n-2}}N$$
as $B_n$-modules, where the action on the right-hand space is 
by left multiplication in $B_n$.
\end{lem}
\begin{proof}
Define a map
$$\phi:B_n\otimes_{B_{n-2}\otimes
  B_2}(N\boxtimes S^{\pm})\rightarrow
B_ne^{\pm}\otimes_{B_{n-2}}N$$
by
$$b\otimes(n\otimes e^{\pm})\longmapsto be^{\pm}\otimes n.$$
We first show that this is well-defined. Let $b=b'b_{n-2}b_2$ for some
$b_{n-2}\in B_{n-2}$ and $b_2\in B_2$. Then
\begin{eqnarray*}
\phi(b\otimes(n\otimes e^{\pm})-b'\otimes(b_{n-2}n\otimes b_2e^{\pm}))&=&
b'b_{n-2}b_2e^{\pm}\otimes n-b'b_2e^{\pm}\otimes b_{n-2}n\\
&=&
b'b_2e^{\pm}b_{n-2}\otimes n-b'b_2e^{\pm}\otimes b_{n-2}n=0
\end{eqnarray*}
as required. The map $\phi$ is clearly a $B_n$-homomorphism. 
We also have a map
$$\psi:B_ne^{\pm}\otimes_{B_{n-2}}N\rightarrow
B_n\otimes_{B_{n-2}\otimes B_2}(N\boxtimes S^{\pm})$$
given by
$$be^{\pm}\otimes n\longmapsto be^{\pm}\otimes (n\otimes e^{\pm}).$$
It is easy to check that $\psi$ is well-defined and that $\psi\phi=\id$
and $\phi\psi=\id$.
\end{proof}

\begin{figure}[ht]
\center{\includegraphics{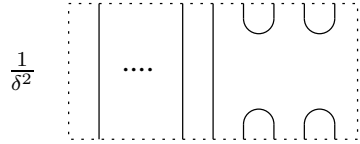}}
\caption{An idempotent in $B_n$}
\label{bige}
\end{figure}

Let $e_{n,4}$ be the idempotent in $B_n$ shown in Figure
\ref{bige}. 

\begin{figure}[ht]
\center{\includegraphics{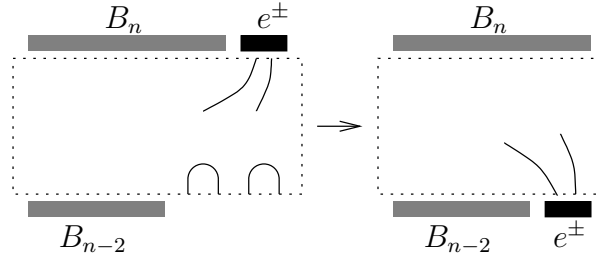}}
\caption{Realising the isomorphism between 
$e^{\pm}B_{n+2}e_{n+2,4}$ and $B_ne^{\pm}$}
\label{e4iso}
\end{figure}

\begin{lem}\label{alt}
 As  left $B_n$- and right $B_{n-2}$-modules we have
$$e^{\pm}B_{n+2}e_{n+2,4}\cong B_ne^{\pm}.$$
\end{lem}
\begin{proof}
Consider the map from $e^{\pm}B_{n+2}e_{n+2,4}$ to
$B_ne^{\pm}$ given on diagrams as shown in Figure
\ref{e4iso}. The grey shaded regions show the actions of $B_n$ from
above, of $B_{n-2}$ from below, and the dark shaded region the
action of the element $e^{\pm}$.
All lines in the diagrams except those indicated remain unchanged;
the two southern arcs in the left-hand diagram are removed, and the
ends of the pair of lines acted on by $e^{\pm}$ are translated
clockwise around the boundary from the northern to the southern side.
This gives an isomorphism of vector spaces, and
clearly preserves the actions of $B_n$ and $B_{n-2}$.
\end{proof}

\begin{cor}\label{pmproj}
The module $\res_n^{\pm}(P_n(\lambda))$ is projective for all
$\lambda\in\Lambda_{n-4}$. 
\end{cor}
\begin{proof}
First note that $B_ne_{n,4}$ is a projective $B_n$-module. Moreover,
as $e_{n,4}B_ne_{n,4}\cong B_{n-4}$ we have that $B_ne_{n,4}$ contains
precisely the indecomposable projectives labelled by elements of
$\Lambda_{n-4}$. By Lemma \ref{alt} we have that
$e^{\pm}B_ne_{n,4}\cong B_{n-2}e^{\pm}$ as left $B_{n-2}$-modules, and
hence $\res_n^{\pm}(P_n(\lambda))$ is projective for all
$\lambda\in\Lambda_{n-4}$.  
\end{proof}

\begin{cor}\label{indrespm}
We have an isomorphism of functors 
$$\ind_n^{\pm}\cong \res_{n+4}^{\pm}G_{n+2}G_{n}.$$
\end{cor}
\begin{proof}
By the definition of $G_n$ and $G_{n+2}$ we have
\begin{eqnarray*}
\res_{n+4}^{\pm}G_{n+2}G_n(N)&=&\res_{n+4}^{\pm}(B_{n+4}e_{n+4,2}
\otimes_{B_{n}} N)\\
&=&e^{\pm}B_{n+4}e_{n+4,4}\otimes_{B_{n}} N\\
&\cong& B_{n+2}e^{\pm}\otimes_{B_{n}} N
\end{eqnarray*}
where the final isomorphism follows from Lemma \ref{alt}. But by Lemma
\ref{altdef} this final module is isomorphic to $\ind_n^{\pm}N$.
\end{proof}

Corollary \ref{indrespm} is an analogue of the relation between
induction, restriction and globalisation in \cite[Lemma 2.6(ii)]{cdm},
corresponding to axiom (A4) for a tower of recollement.

Given two partitions $\lambda$ and $\mu$, we write
$\lambda\rhd\rhd^+\mu$, or $\mu\lhd\lhd^+\lambda$, if $\mu$ can be
obtained from $\lambda$ by removing two boxes and $\lambda/\mu$ is not
the partition $(1,1)$. Similarly we write $\lambda\rhd\rhd^-\mu$, or
$\mu\lhd\lhd^-\lambda$  if $\mu$ can be
obtained from $\lambda$ by removing two boxes and $\lambda/\mu\neq
(2)$. We will write $\lambda\lhd\rhd\mu$ if $\mu$ is obtained from
$\lambda$ by removing a box and then adding a box.

The next theorem describes the structure of
$\res_n^{\pm}\Delta_n(\lambda)$, and so is an analogue of the usual
induction and restriction rules in \cite[Theorem 4.1 and Corollary
  6.4]{dhw} (and use the same strategy for the proof).

\begin{thm}\label{bigres}
Suppose that $\lambda$ is a partition of $m=n-2t$ for some $t\geq 0$.\\
(i) There is a
filtration of $B_{n-2}$-modules
$$W_0\subseteq W_1\subseteq W_2=\res_{n}^{\pm}\Delta_n(\lambda)$$
with
$$W_0\cong \bigoplus_{\mu\lhd\lhd^{\pm}\lambda}\Delta_{n-2}(\lambda)
\quad\quad
W_1/W_0\cong \bigoplus_{\mu\lhd\rhd\lambda}\Delta_{n-2}(\lambda)
\quad\quad
W_2/W_1\cong
\bigoplus_{\mu\rhd\rhd^{\pm}\lambda}\Delta_{n-2}(\lambda)$$
where any $\Delta_{n-2}(\mu)$ which does not make sense is taken as
$0$.\\
(ii) There is a 
filtration of $B_{n+2}$-modules
$$U_0\subseteq U_1\subseteq U_2=\ind_{n}^{\pm}\Delta_n(\lambda)$$
with
$$U_0\cong \bigoplus_{\mu\lhd\lhd^{\pm}\lambda}\Delta_{n+2}(\lambda)
\quad\quad
U_1/U_0\cong \bigoplus_{\mu\lhd\rhd\lambda}\Delta_{n+2}(\lambda)
\quad\quad
U_2/U_1\cong
\bigoplus_{\mu\rhd\rhd^{\pm}\lambda}\Delta_{n+2}(\lambda).$$
\end{thm}
\begin{proof}
Part (ii) follows from part (i) by Corollary \ref{indrespm}. For the
rest of the proof we will work with the concrete realisation of
standard modules given in Section \ref{basics}.  By definition we have
$$\res_n^{\pm}\Delta_n(\lambda)=e^{\pm}I_n^t\otimes_{\Sigma_m}S^{\lambda}$$
and we will represent an element $e^{\pm}X_{w,1,id}\otimes x$ in this
space diagrammatically as shown in Figure \ref{howrep}.

\begin{figure}[ht]
\center{\includegraphics{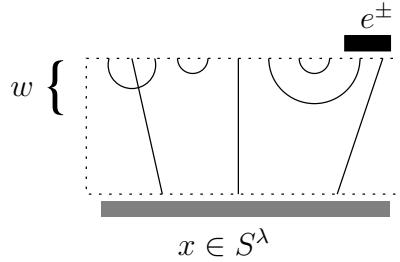}}
\caption{Representing the element 
$e^{\pm}X_{w,1,id}\otimes x$ in $e^{\pm}I_n^t\otimes_{\Sigma_m} S^{\lambda}$}
\label{howrep}
\end{figure}

We are now in a position to define the various spaces $W_0$, $W_1$,
and $W_2$. Choose a fixed basis $V(\lambda)$ for $S^{\lambda}$ and set 
$$V_{n,t}^0=\{w\in V_{n,t}:\mbox{\rm $n-1$ and $n$ are free in
  $w$}\}$$
$$V_{n,t}^1=\{w\in V_{n,t}:\mbox{\rm $n-1$ is on an arc and $n$
 is free in $w$}\}$$
$$V_{n,t}^2=\{w\in V_{n,t}:\mbox{\rm $n-1$ is linked to $j$ and
$n$ is linked to $i$ in $w$ with $i<j\leq n-2$}\}.$$
Then for $0\leq i\leq 2$ we set
$$W_i=\spann\{e^{\pm}X_{w,1,id}\otimes x: w\in V_{n,t}^j\  \mbox{\rm with }
j\leq i\ \mbox{\rm and } x\in V(\lambda)\}$$

Note that if $w\in V_{n,t}^1\cup V_{n,t}^2$ and $w'$ is obtained from
$w$ by swapping nodes $n-1$ and $n$, then
$$e^{\pm}X_{w,1,id}=\pm e^{\pm}X_{w',1,id}.$$
Moreover, if there is an arc linking nodes $n-1$ and $n$ in $w$ then 
$$e^{\pm}X_{w,1,id}=0.$$
Thus we have that $$W_2=e^{\pm}\Delta_n(\lambda)$$
and $W_0$ and $W_1$ are submodules of $e^{\pm}\Delta_n(\lambda)$.

We first show that
\begin{equation}\label{firststep}
W_0\cong I_{n-2}^t\otimes_{\Sigma_{m-2}}\sigma^{\pm}S^{\lambda}
\end{equation}
where $\sigma^{\pm}$ represents the symmetriser/antisymmetriser on the
last two lines in $\Sigma_m$ and $\Sigma_{m-2}\subset \Sigma_m$ acts on
the first $m-2$ lines. Note that
\begin{eqnarray*}
\sigma^{\pm}S^{\lambda}
&=&\sigma^{\pm}\res_{\Sigma_{m-2}\times\Sigma_2}^{\Sigma_m}S^{\lambda}\\
&=&\sigma^{\pm}\left(\bigoplus_{\mu\vdash m-2,\ \nu\vdash 2}
c_{\mu,\nu}^{\lambda}(S^{\mu}\boxtimes S^{\nu})\right)
=\bigoplus_{\mu\vdash m-2}c_{\mu,*}^{\lambda}S^{\mu} 
\end{eqnarray*}
where $*$ equals $(2)$ for $\sigma^+$ and $(1,1)$ for $\sigma^-$. As 
$$
c_{\mu,(2)}^{\lambda}=\left\{\begin{array}{ll}
1 & \wif \mu\lhd\lhd^+\lambda\\
0&\otherwise\end{array}\right.\quad\quad\wand 
\quad\quad
c_{\mu,(1,1)}^{\lambda}=\left\{\begin{array}{ll}
1 & \wif \mu\lhd\lhd^-\lambda\\
0&\otherwise\end{array}\right.$$
it will follow from (\ref{concrete}) that
$$W_0\cong \bigoplus_{\mu\lhd\lhd^{\pm}\lambda}\Delta_{n-2}(\lambda)$$
as required.

Note that for $w\in V_{n,t}^0$ the lines from $n-1$ and $n$ are
propagating in $X_{w,1,id}$, and so we have
$$e^{\pm}X_{w,1,id}\otimes S^{\lambda}=
X_{w,1,id}e^{\pm}\otimes S^{\lambda}
= X_{w,1,id}\otimes \sigma^{\pm}S^{\lambda}.$$
For $w\in V_{n,t}^0$ define $\overline{w}\in V_{n-2,t}$ by removing nodes
$n-1$ and $n$, and a map 
$$\phi_0: W_0\rightarrow
I_{n-2}^t\otimes_{\Sigma_{m-2}}\sigma^{\pm}S^{\lambda}$$
by
$$e^{\pm}X_{w,1,id}\otimes x=X_{w,1,id}\otimes \sigma^{\pm}x
\longmapsto X_{\overline{w},1,id}\otimes\sigma^{\pm}x.$$ 
It is clear that $\phi_0$ is an isomorphism of vector spaces, and
commutes with the action of $B_{n-2}$. This proves (\ref{firststep}).

Next we will show that
\begin{equation}\label{secondstep}
W_1/W_0\cong I_{n-2}^{t-1}\otimes_{\Sigma_{m}}
\ind_{\Sigma_{m-1}}^{\Sigma_m}\res_{\Sigma_{m-1}}^{\Sigma_m}S^{\lambda}.
\end{equation}
Note that
$$\ind_{\Sigma_{m-1}}^{\Sigma_m}\res_{\Sigma_{m-1}}^{\Sigma_m}S^{\lambda}
=\ind_{\Sigma_{m-1}}^{\Sigma_m}
\bigl(\bigoplus_{\nu\lhd\lambda}
S^{\nu}\bigr)=\bigoplus_{\nu\lhd\lambda}
\bigl(\ind_{\Sigma_{m-1}}^{\Sigma_m}S^{\nu}\bigr)=
\bigoplus_{\mu\lhd\rhd\lambda}S^{\mu}
$$
and so it will follow from (\ref{concrete}) that
$$
W_1/W_2\cong \bigoplus_{\mu\lhd\rhd\lambda}\Delta_{n-2}(\lambda)
$$
as required.

We will need an explicit description of
$\ind_{\Sigma_{m-1}}^{\Sigma_m}\res_{\Sigma_{m-1}}^{\Sigma_m}S^{\lambda}$. The
quotient $\Sigma_m/\Sigma_{m-1}$ has coset representatives
$$\{\tau_i=(i,m):1\leq i\leq m\}$$
where $(m,m)=1$. Therefore
$\ind_{\Sigma_{m-1}}^{\Sigma_m}\res_{\Sigma_{m-1}}^{\Sigma_m}S^{\lambda}$
has a basis 
$$\{(i,x): 1\leq i\leq m, x\in V(\lambda)\}$$
and the action of $\theta\in\Sigma_m$ is given by
$$\theta(i,x)=(j,\theta'x)$$ where $\theta\tau_i=\tau_j\theta'$ for a
unique $1\leq j\leq m$ and $\theta'\in\Sigma_{m-1}$.

For $1\leq i\leq m$ set 
$$\sigma_i=(i,m,m-1,m-2,\ldots,i+1)$$ and
for $w\in V^1_{n,t}$ define $\overline{w}\in V_{n-2.t-1}$ by removing the
nodes $n-1$ and $n$ and removing the arc from $n-1$ (which will thus
introduce a new free node elsewhere in $\overline{w}$). Now we can define a
map
$$\phi_1: W_1/W_0\rightarrow I_{n-2}^{t-1}\otimes_{\Sigma_{m}}
\ind_{\Sigma_{m-1}}^{\Sigma_m}\res_{\Sigma_{m-1}}^{\Sigma_m}S^{\lambda}$$
by
$$e^{\pm}X_{w,1,id}\otimes x
\longmapsto X_{\overline{w},1,id}\sigma_i\otimes(m,x)$$
if node $n-1$ is linked to node $i$ in $w$. This is illustrated
graphically in Figure \ref{phi1}.

\begin{figure}[ht]
\center{\includegraphics{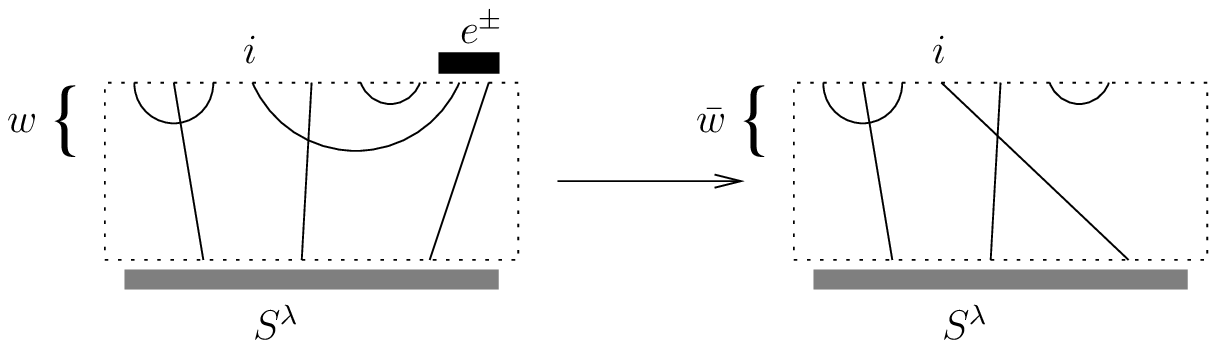}}
\caption{An example of the effect of the map $\phi_1$}
\label{phi1}
\end{figure}

Note that for every $v\in V_{n-2,t-1}$ there are exactly $m$ elements
$w\in V^1_{n,t}$ satisfying $\overline{w}=v$, as $n-1$ can be joined to any
of the $m$ free vertices in $v$. Note also that
$$\sigma_i=(i,m)(m-1,m-2,\ldots,i+1,i)=(i,m)\sigma_i'$$
where $\sigma_i'\in\Sigma_{m-1}$, and so $\sigma_i(m,x)=(i,\sigma_i'x)$.

Given $v\in V_{n-2,t-1}$, $1\leq i\leq m$, and $x\in V(\lambda)$ pick
$w\in V^1_{n,t}$ with $\overline{w}=v$ and $n-1$ joined to the $i$th free
node. Then
\begin{eqnarray*}
\phi_1(e^{\pm}X_{w,1,id}\otimes (\sigma_i')^{-1}x)&=&
X_{v,1,id}\sigma_i\otimes (m,(\sigma_i')^{-1}x)\\
&=& X_{v,1,id}\otimes (i,\sigma_i'(\sigma_i')^{-1}x)
= X_{v,1,id}\otimes (i,x)
\end{eqnarray*}
and so $\phi_1$ is surjective. Moreover
$$\dim W_1/W_2=m|V_{n-2,t-1}|\dim S^{\lambda}=\dim I_{n-2}^{t-1}\otimes
\ind_{\Sigma_{m-1}}^{\Sigma_m}\res_{\Sigma_{m-1}}^{\Sigma_m}S^{\lambda}
$$
and so $\phi_1$ is an isomorphism of vector spaces. It remains to show
that $\phi_1$ commutes with the action of $B_{n-2}$.

First consider the action of $\tau\in \Sigma_{n-2}$. The actions of
$\phi_1$ and $\tau$ can be seen to commute by the schematic diagram in
Figure \ref{comm1}, noting that $\tau(\overline{w})=\overline{\tau(w)}$.

\begin{figure}[ht]
\center{\includegraphics{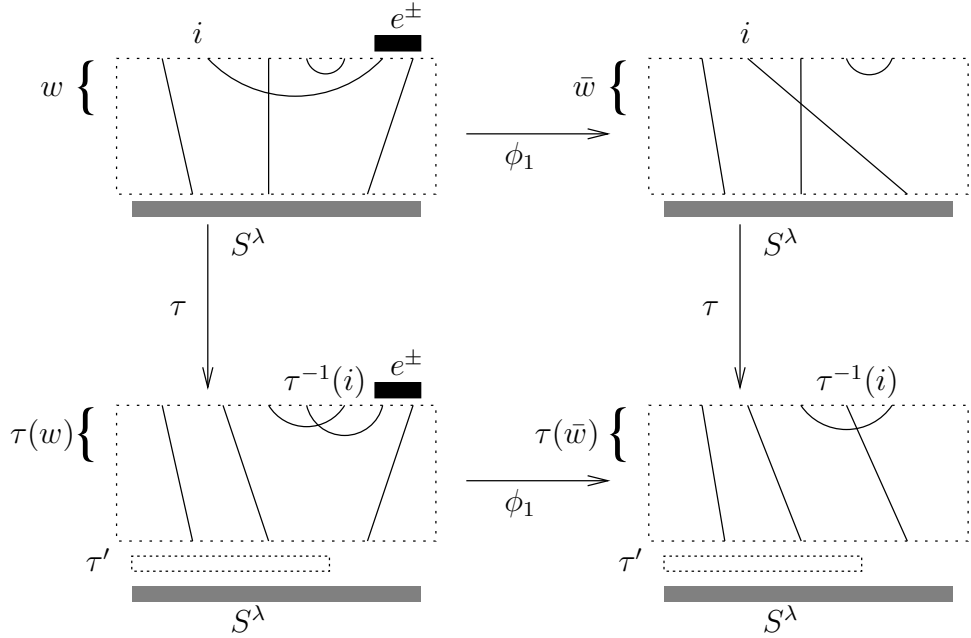}}
\caption{A diagrammatic illustration that $\phi_1\tau=\tau\phi_1$}
\label{comm1}
\end{figure}

Next consider the action of $X_{jk}\in B_{n-2}$. If $j,k\neq i$ then
it is clear that $X_{jk}$ commutes with $\phi_1$. Now consider the
action of $X_{ij}$. There are two cases: (i) $j$ is a free node in $w$,
and (ii) $j$ is linked to some node $k$ in $w$.  Case (i) is
illustrated schematically in Figure \ref{comm2}. The lower left
diagram in Figure \ref{comm2} represents $0$ as it lies in $W_0$. The
lower right diagram represents $0$ as there is a decrease in the
number of propagating lines. Therefore the dotted arrow is an equality
and the diagram commutes.

\begin{figure}[ht]
\center{\includegraphics{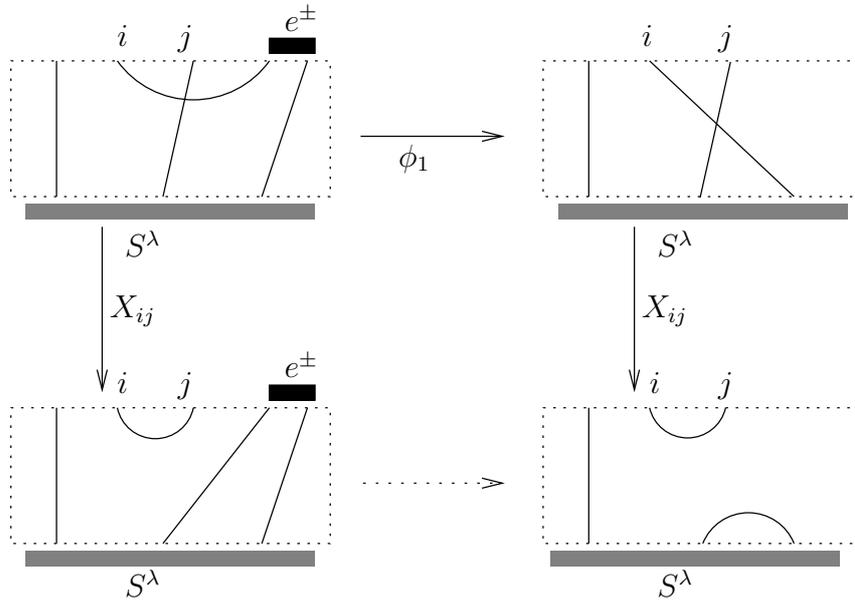}}
\caption{The action of $X_{ij}$ and $\phi_1$: case (i)}
\label{comm2}
\end{figure}

Case (ii) is illustrated schematically in Figure \ref{comm3}. Again we
see that $X_{ij}$ commutes with the action of $\phi_1$, and so we have
shown that $\phi_1$ is a $B_{n-2}$-isomorphism. This completes the
proof of (\ref{secondstep}).

\begin{figure}[ht]
\center{\includegraphics{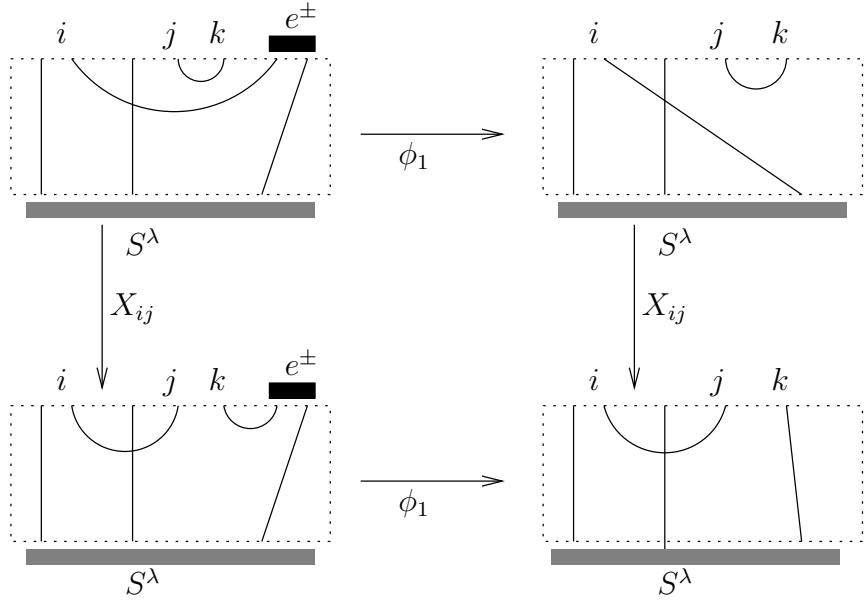}}
\caption{The action of $X_{ij}$ and $\phi_1$: case (ii)}
\label{comm3}
\end{figure}

Finally, we shall show that 
\begin{equation}\label{thirdstep}
W_2/W_1\cong I_{n-2}^{t-2}\otimes_{\Sigma_{m-2}}
\ind_{\Sigma_m\times \Sigma_2}^{\Sigma_{m+2}}(S^{\lambda}\boxtimes
S^{\pm}).
\end{equation}
As 
$$\ind_{\Sigma_m\times
  \Sigma_2}^{\Sigma_{m+2}}(S^{\lambda}\boxtimes S^{\pm})\cong
\bigoplus_{\mu\rhd\rhd^{\pm}\lambda}S^{\mu}$$
it will follow from (\ref{concrete}) that
$$W_2/W_1\cong \bigoplus_{\mu\rhd\rhd^{\pm}\lambda}\Delta_{n-2}(\lambda)$$
which will complete the proof. 

We will need an explicit description of
$\ind_{\Sigma_{m}\times \Sigma_2}^{\Sigma_{m+2}}(S^{\lambda}\boxtimes
S^{\pm})$. The
quotient $\Sigma_{m+2}/(\Sigma_{m}\times\Sigma_2)$ has coset representatives
$$\{\tau_{ij}=(i,m+1)(j,m+2):1\leq i<j\leq m+2\}$$
where $(m+1,m+1)=(m+2,m+2)=1$. Therefore
$\ind_{\Sigma_{m}\times\Sigma_2}^{\Sigma_{m+2}}(S^{\lambda}\boxtimes S^{\pm})$
has a  basis 
$$\{(i,j;x\otimes\sigma^{\pm}): 1\leq i<j\leq m+2, x\in V(\lambda)\}$$
and the action of $\theta\in\Sigma_{m+2}$ is given by
$$\theta(i,j;x\otimes\sigma^{\pm})=(k,l;\theta'(x\otimes\sigma^{\pm}))$$ 
where $\theta\tau_{ij}=\tau_{kl}\theta'$ for a
unique $1\leq k<l\leq m+2$ and $\theta'\in\Sigma_{m}\times\Sigma_2$.

For $1\leq r<s\leq m+2$ set
$$\sigma_{r,s}=(r,s,s-1,s-2,\ldots,r+1)$$ and for $w\in V^2_{n,t}$
define $\overline{w}\in V_{n-2,t-2}$ by removing the nodes $n-1$ and
$n$ and removing the arcs from $n-1$ and $n$ (which will thus
introduce two new free nodes elsewhere in $\overline{w}$). Now we can
define a map
$$\phi_2: W_2/W_1\rightarrow I_{n-2}^{t-2}\otimes_{\Sigma_{m+2}}
\ind_{\Sigma_m\times\Sigma_2}^{\Sigma_{m+2}}(S^{\lambda}\boxtimes S^{\pm})$$
by
$$e^{\pm}X_{w,1,id}\otimes x \longmapsto
X_{\overline{w},1,id}\sigma_{j,m+2}\sigma_{i,m+1}
\otimes(m+1,m+2;,x\otimes\sigma^{\pm})$$ if $n-1$ is linked to $j$ and
$n$ is linked to $i$ in $w$. This is illustrated graphically in Figure
\ref{phi2}.

\begin{figure}[ht]
\center{\includegraphics{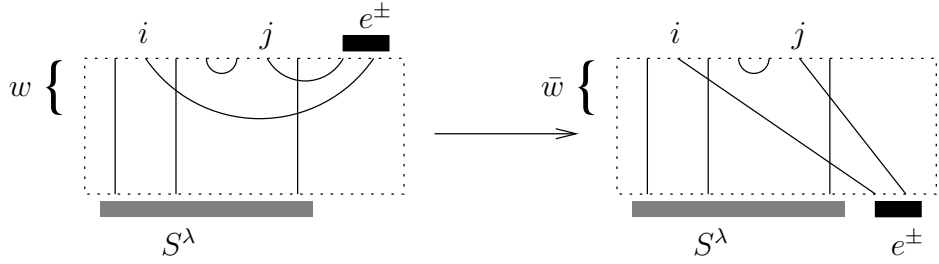}}
\caption{An example of the effect of the map $\phi_2$}
\label{phi2}
\end{figure}

Arguing as for $\phi_1$ we can show that $\phi_2$ is an isomorphism of
vector spaces. Thus we will be done if we can show that $\phi_2$
commutes with the action of $B_{n-2}$.

First consider the action of $\tau\in\Sigma_{n-2}$. The actions of
$\phi_2$ and $\tau$ are illustrated schematically in Figure
\ref{comm4}. Again we use the fact that
$\tau(\overline{m})=\overline{\tau(w)}$, while in the bottom pair of
diagrams we have used the action of $e^{\pm}$ on each side, which in
each case gives a coefficient of $\pm1$. We see that the actions of
$\phi_2$ does commute with $\tau$ as required.

\begin{figure}[ht]
\center{\includegraphics{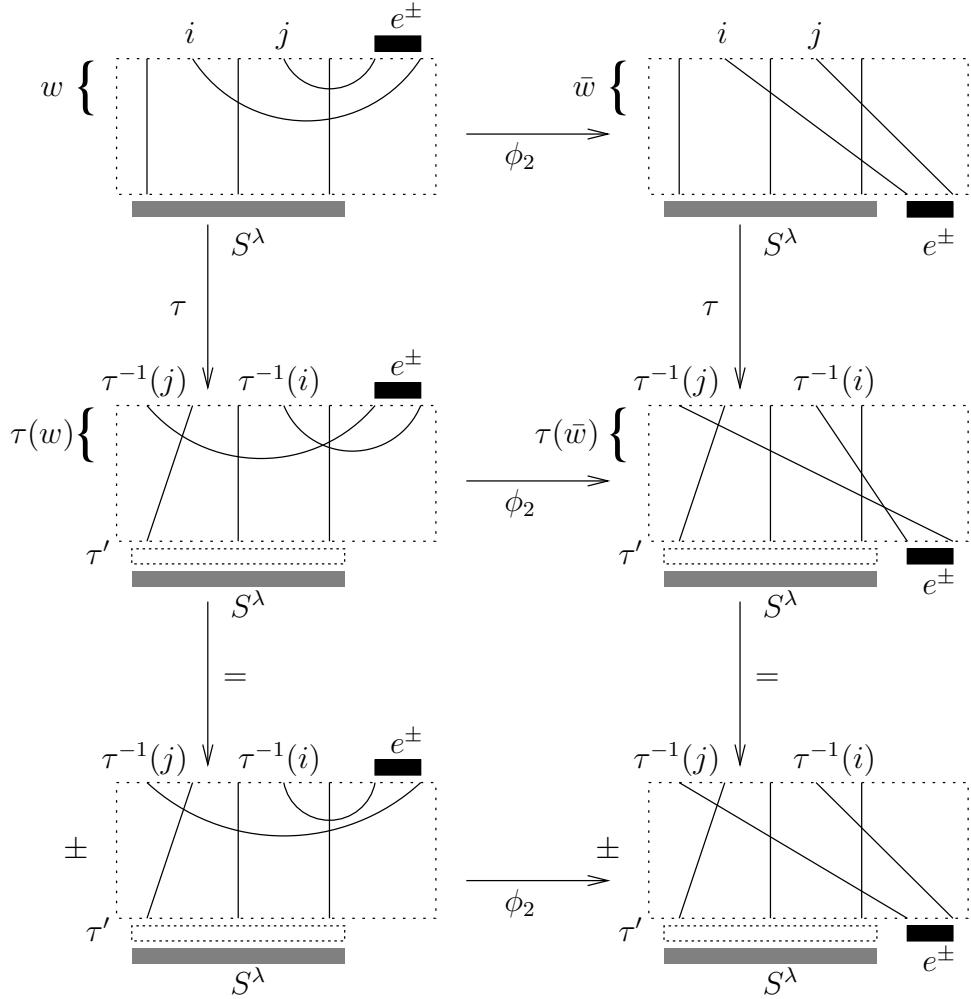}}
\caption{A diagrammatic illustration that $\phi_2\tau=\tau\phi_1$}
\label{comm4}
\end{figure}

It remains to check that $\phi_2$ commutes with the action of
$X_{kl}\in B_{n-2}$. If $\{k,l\}$ is disjoint from $\{i,j\}$ then it
is clear that $X_{k,l}$ commutes with $\phi_2$. If $k=i$ and $l=j$ it
is easy to verify that
$$X_{ij}e^{\pm}X_{w,1,id}=0\quad\wand\quad
X_{i,j}\phi_2(e^{\pm}X_{w,1,id}\otimes x)=0.$$
Thus we just have to check what happens when $k=i$ and $l\neq
j$. There are two cases: (i) $l$ is a free node in $w$, and (ii) $l$
is linked to some node $h$ in $w$.

Case (i) is illustrated schematically in Figure \ref{comm5}. The lower left
diagram in Figure \ref{comm5} represents $0$ as it lies in $W_1$. The
lower right diagram represents $0$ as there is a decrease in the
number of propagating lines. Therefore the dotted arrow is an equality
and the diagram commutes.

\begin{figure}[ht]
\center{\includegraphics{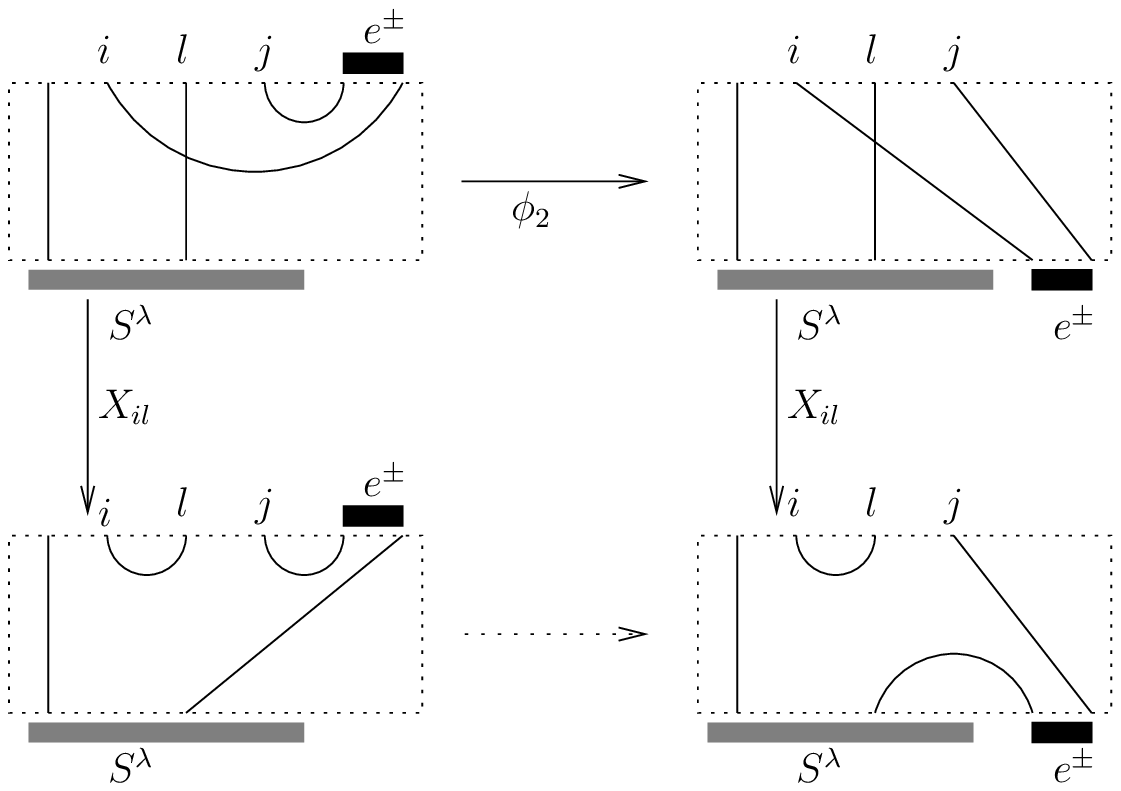}}
\caption{The action of $X_{il}$ and $\phi_2$: case (i)}
\label{comm5}
\end{figure}

\begin{figure}[ht]
\center{\includegraphics{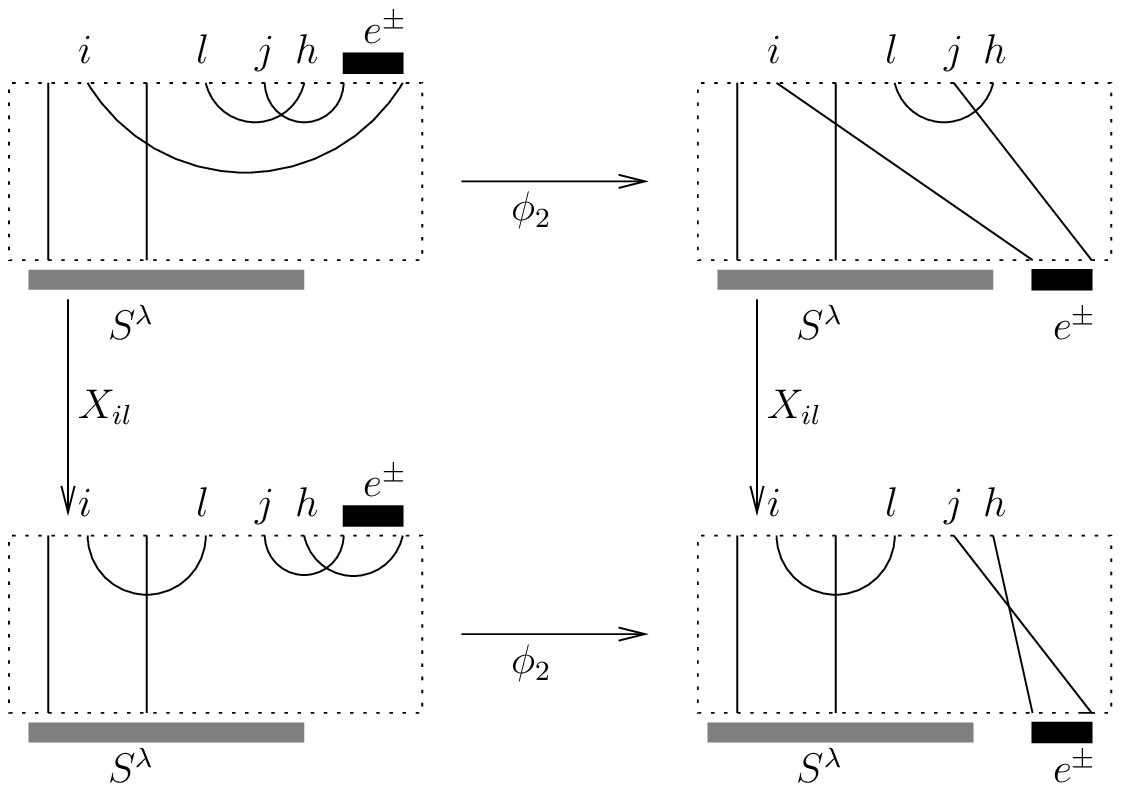}}
\caption{The action of $X_{il}$ and $\phi_2$: case (ii)}
\label{comm6}
\end{figure}

Case (ii) is illustrated schematically in Figure \ref{comm6}. Again we
see that $X_{il}$ commutes with the action of $\phi_2$, and so we are done.
\end{proof}

\section{Translation equivalence for the Brauer algebra}\label{brauertran}

In Section \ref{towermor} we saw how translation equivalence of
weights for a tower of recollement implies Morita equivalences of the
corresponding blocks (when the blocks are truncated to contain the
same number of simples).  We will now reinterpret this in the language
of alcove geometry in the case of the Brauer algebra.

Given a partition $\lambda$, we saw in (\ref{resrule}) that the set
$\supp(\lambda)$ consists of those partitions obtained from $\lambda$
by the addition or subtraction of a box from $\lambda$, all with
multiplicity one.  We denote by $\lambda\pm\epsilon_i$ the composition
obtained by adding/subtracting a box from row $i$ of $\lambda$.

\begin{lem}\label{nosep}
If $\lambda'=\lambda\pm\epsilon_i$ then there cannot exist a
reflection hyperplane separating $\lambda$ from $\lambda'$.
\end{lem}
\begin{proof}
Suppose that $R$ is a reflection hyperplane between $\lambda$ and
$\lambda'$, and denote their respective reflections by $r(\lambda)$
and $r(\lambda')$. Either $\lambda$ is the reflection of $\lambda'$ or
the line from $\lambda$ to $\lambda'$ is not orthogonal to the
hyperplane.

The former case is impossible as two weights differing by
one box cannot be in the same block. For the latter case, note that
the distance between $\lambda$ and $\lambda'$ is one. Therefore at
least one of the distances from $\lambda$ to $r(\lambda)$ and
$\lambda'$ to $r(\lambda')$ is less than one. But this is impossible,
as $r(\lambda)$ and $r(\lambda')$ are also elements of the lattice of
weights. 
\end{proof}

Given a facet $F$, we denote by $\overline{F}$ the closure of $F$ in
$\EE^{\infty}$. This will consist of a union of facets.

\begin{lem}\label{onwall}
If $\lambda'\in\supp(\lambda)$ and $\lambda'\in F$ for some facet $F$
then
$$|{\mathcal B}(\lambda')\cap\supp(\lambda)|>1$$ 
if and only if $\lambda\in\overline{F}\backslash F$.
\end{lem}
\begin{proof}
We first show that if $|{\mathcal B}(\lambda')\cap\supp(\lambda)|>1$
then $\lambda\in\overline{F}\backslash F$.  By the interpretation of
blocks in terms of contents of partitions (Theorem \ref{balver}) there
is precisely one other weight $\lambda''$ in $\supp(\lambda)$ in the
same block as $\lambda'$. Also, one of these weights is obtained from
$\lambda$ by adding a box, and one by subtracting a box.  But this
implies that $\lambda'$ is the reflection of $\lambda''$ about some
hyperplane $H$; this reflection must fix the midpoint on the line from
$\lambda'$ to $\lambda''$, which is $\lambda$, and so $\lambda\in H$.

By Lemma \ref{nosep} there is no hyperplane separating $\lambda$ from
$\lambda'$. However to complete the first part of the proof we still
need to show that if $\lambda'\in H'$ for some hyperplane $H'$ then
$\lambda\in H'$ too. First note that by Proposition \ref{Wis}(ii) $H$
and $H'$ must be $(i,j)_-$ and $(k,l)_-$ hyperplanes respectively for
some quadruple $i,j,k,l$. It is easy to check that either $(i,j)_-$
fixes $H'$ or $(i,j)_-\cdot_{\delta} H'$ is an
$(i,l)$-hyperplane. But $\lambda''\in(i,j)_-\cdot_{\delta}H'$ is
dominant and so $(i,j)_-$ must fix $H'$. Hence $\lambda''\in H'$ and
as $\lambda$ is the midpoint between $\lambda'$ and $\lambda''$ we
must have $\lambda\in H'$.

For the reverse implication, suppose that $\lambda\in
\overline{F}\backslash F$. Then for all hyperplanes $H'$ with
$\lambda'\in H'$ we have $\lambda\in H'$ and there is (at least) one
hyperplane $H$ with $\lambda\in H$ and $\lambda'\notin H$. Suppose
that $H$ is an $(i,j)_-$-hyperplane, and consider
$\lambda''=(i,j)_-\cdot_{\delta}\lambda'$. 

If $\lambda''\in X^+$ then we are done. Otherwise by Lemma \ref{nosep}
we have that $\lambda''$ must lie on the boundary of the dominant
region, and hence in some $(k,l)$-hyperplane $\tilde{H}$. Now
$(i,j)_-\cdot_{\delta}\lambda''=\lambda'$ and hence
$\lambda'\in(i,j)_-\cdot_{\delta}\tilde{H}=H'\neq \tilde{H}$ (as
$\lambda'\in X^+$). Therefore we must have $\lambda\in H'$. But
$\lambda$ is fixed by $(i,j)_-$ and so $\lambda\in\tilde{H}\cap
H'$. This implies that $\lambda\notin X^+$ which is a contradiction.
Thus we have shown that $\lambda''\in X^+$ and so $|{\mathcal
  B}(\lambda')\cap\supp(\lambda)|>1$.
\end{proof}

\begin{thm}\label{alcovemor}
If $\lambda$ is in an alcove then $\mu$ is in the same translation
class as $\lambda$ if and only if it is in the same alcove.
\end{thm}
\begin{proof} 
By Lemmas \ref{nosep} and \ref{onwall} it is enough to show that 
if $\mu$ is in the same alcove as $\lambda$ then $\mu$ can be obtained
from $\lambda$ by repeatedly adding or subtracting a box without ever
leaving this alcove.

Suppose that $\lambda$ and $\mu$ are in the same alcove, and set
$x=\lambda+\rho_{\delta}$ and $y=\mu+\rho_{\delta}$, the corresponding
vectors in $A^+$.  Recall that there is a permutation $\pi$ defining
the alcove $A$ introduced in Section \ref{weyl}. We may assume that
$|x_{\pi(1)}|\leq |y_{\pi(1)}|$. Consider the sequence obtained by
repeatedly adding (or subtracting) $1$ from $y_{\pi(1)}$ until we
obtain $x_{\pi(1)}$. At each stage the vector $v$ obtained is of the
form $\tau+\rho_{\delta}$ for some weight $\tau$, and the sequence of weights
thus obtained are such that each consecutive pair are translation
equivalent.  Now we repeat the process to convert $y_{\pi(2)}$ into
$x_{\pi(2)}$ (note that $y_{\pi(2)}$ and $x_{\pi(2)}$ have the same
sign, and so the chain of weights constructed will always have
$\pi(2)$-coordinate satisfying the defining conditions for the
alcove). We continue in this manner until we have converted $y$ into
$x$. This constructs a chain of translation equivalent weights
connecting $\lambda$ and $\mu$ and so we are done.
\end{proof}

\begin{rem} Theorem \ref{alcovemor} shows that the geometry on the
  weight space for $B_n(\delta)$ comes naturally from the induction
  and restriction functors when the alcoves are non-empty (i.e. for
  $\delta>0$).
\end{rem}

We would like to extend Theorem \ref{alcovemor} to the case of two
weights in the same facet. However, not all weights in the same facet
are in the same translation class. To see this, note that a hyperplane
is defined by the equation $x_i=-x_j$ for some fixed pair $i$ and
$j$. Any modification of a weight in such a hyperplane by
adding or subtracting a single box cannot alter the value of the $i$th
or $j$th coordinate without leaving the hyperplane. However, we will
see that if we also use the modified translation functors introduced
in Section \ref{genir} then we do get the desired equivalences within
facets.

Let $\supp^2(\lambda)=\supp(\supp(\lambda))$. This set consists of
those partitions obtained from $\lambda$ by adding two boxes, removing
two boxes, or adding a box and removing a box.

\begin{lem}\label{wallcase}
Suppose that $\lambda, \tilde{\lambda}\in X^+$ with
$\tilde{\lambda}\in\supp^2(\lambda)$. If $\lambda$ and
$\tilde{\lambda}$ are in the same facet then  
$$\supp^2(\lambda)\cap{\mathcal B}(\tilde{\lambda})=\{\tilde{\lambda}\}.$$
\end{lem}
\begin{proof}
We take $\lambda'\in\supp^2(\lambda)$ with $\lambda'\neq
\tilde{\lambda}$ and show that the above assumptions imply that
$\lambda'\notin{\mathcal B}(\lambda)$. There are six possible cases.

(i) Suppose that $\tilde{\lambda}=\lambda-\epsilon_i+\epsilon_j$ and
$\lambda'=\lambda-\epsilon_k+\epsilon_l$. For these two weights to be
in the same block the boxes $\epsilon_i$ and $\epsilon_l$ must pair up
(and so must $\epsilon_j$ and $\epsilon_k$) in the sense of condition
(1) for a balanced partition. This implies that there is a simple
reflection $(i,l)_-$ taking $\lambda-\epsilon_i$ to
$\lambda+\epsilon_l$, which fixes $\lambda$. Hence $\lambda$ is in the
$(i,l)_-$-hyperplane. However, $\tilde{\lambda}$ is not in this
hyperplane, contradicting our assumption that they are in the same facet.

(ii) Suppose that $\tilde{\lambda}=\lambda+\epsilon_i+\epsilon_j$ and
$\lambda'=\lambda+\epsilon_k+\epsilon_l$. For these two weights to be
in the same block the elements
$\epsilon_i,\epsilon_j,\epsilon_k,\epsilon_l$ must all be distinct and
$\epsilon_i$ and $\epsilon_j$ must pair up (and so must $\epsilon_k$
and $\epsilon_l$). Thus there is a reflection taking $\lambda$ to
$\tilde{\lambda}$, which contradicts our assumption.

(iii) Suppose that $\tilde{\lambda}=\lambda-\epsilon_i-\epsilon_j$ and
$\lambda'=\lambda-\epsilon_k-\epsilon_l$. This is similar to (ii).

(iv) Suppose that $\tilde{\lambda}=\lambda-\epsilon_i+\epsilon_j$ and
$\lambda'=\lambda+\epsilon_k+\epsilon_l$. For these two weights to be
in  the same block we must have $j=l$ (say). But then
$\lambda+\epsilon_k+\epsilon_j$ is the reflection of
$\lambda-\epsilon_i+\epsilon_j$ through the $(i,k)_-$-hyperplane, and
hence $\lambda+\epsilon_j$ is in the $(i,k)_-$-hyperplane. Therefore
$\lambda$ is also in this hyperplane, but $\tilde{\lambda}$ is not,
which gives a contradiction.

(v) Suppose that $\tilde{\lambda}=\lambda-\epsilon_i+\epsilon_j$ and
$\lambda'=\lambda-\epsilon_k-\epsilon_l$. This is similar to (iv).

(vi) Suppose that $\tilde{\lambda}=\lambda+\epsilon_i+\epsilon_j$ and
$\lambda'=\lambda-\epsilon_k-\epsilon_l$. First note that $\epsilon_i$
and $\epsilon_j$ cannot pair up (as this would imply that $\lambda$
and $\tilde{\lambda}$ are not in the same facet). So for these two
weights to be in the same block we must have $\epsilon_i$ pairing up
with $\epsilon_k$ (say) and $\epsilon_j$ pairing up with
$\epsilon_l$. But then $\lambda-\epsilon_k$ is the reflection of
$\lambda+\epsilon_i$ through the $(i,k)_-$-hyperplane, which implies
that $\lambda$ is in this hyperplane but $\tilde{\lambda}$ is not,
which gives a contradiction.
\end{proof}

\begin{lem}\label{samen}
Suppose that $\lambda,\tilde{\lambda}\in X^+$ with
$\tilde{\lambda}=\lambda-\epsilon_i+\epsilon_j\in\supp^2(\lambda)$,
and that $\lambda$ lies on the $(ij)_-$-hyperplane. Then $\lambda$ and
$\tilde{\lambda}$ are in the same facet. Moreover, if
$\mu=w\cdot_{\delta}\lambda\in X^+$ for some $w\in W$ then
$\tilde{\mu}=w\cdot_{\delta}\tilde{\lambda}$ satisfies
$$\tilde{\mu}=\mu-\epsilon_s+\epsilon_t$$
for some $s,t$. 
\end{lem}
\begin{proof}
The fact that $\lambda$ and $\tilde{\lambda}$ are in the same facet is
clear. Now suppose that $\mu=w\cdot_{\delta}\lambda$ and
$\tilde{\mu}=w\cdot_{\delta}\tilde{\lambda}$. Then
$\tilde{\mu}=\mu+\beta$ where $\beta=\pm(\epsilon_s+\epsilon_t)$ or
$\beta=\epsilon_s-\epsilon_t$ for some $s,t$. Suppose for a
contradiction that $\beta=\pm(\epsilon_s+\epsilon_t)$. Note that $\mu$
and $\tilde{\mu}$ are in the same facet, and so for any
$(k,l)_-$-hyperplane on which $\mu$ and $\tilde{\mu}$ lie we must have
that $\beta=\tilde{\mu}-\mu$ lies on the \emph{unshifted}
$(k,l)_-$-hyperplane. This implies that $s\neq k,l$ and $t\neq k,l$.

We have a sequence of dominant weights $\mu$, $\mu'=\mu\pm\epsilon_s$
and $\tilde{\mu}=\mu\pm(\epsilon_s+\epsilon_t)$ which are each at distance
$1$ from their neighbours in the sequence. We have already seen that
they all lie on the same set of hyperplanes. Moreover, by Lemma
\ref{nosep} there cannot exist a hyperplane separating $\mu$ from
$\mu'$ or $\mu'$ from $\tilde{\mu}$. So $\mu$, $\mu'$ and
$\tilde{\mu}$ all lie in the same facet.

Now consider the image of these three weights under $w^{-1}$. We get a
corresponding sequence $\lambda$, $\lambda'$ and
$\tilde{\lambda}$. These weights must also lie in a common facet (and
hence $\lambda'$ is dominant) and are distance $1$ from their
neighbours. This forces $\lambda'=\lambda-\epsilon_i$ or
$\lambda'=\lambda+\epsilon_j$. However $\lambda'$ cannot be in the
same facet as $\lambda$ as it does not lie on the
$(i,j)_-$-hyperplane, which gives the desired contradiction.
\end{proof}

Let $\res_n^{\lambda,\pm}=\pr^{\lambda}\res_n^\pm$ and
$\ind_n^{\lambda,\pm}=\pr^{\lambda}\ind_n^{\pm}$.
We say that $\lambda$ and $\mu$ are in the same \emph{$(\pm)$-translation
class} if there is a chain of dominant weights 
$$\lambda=\lambda^0,\lambda^1,\ldots,\lambda^r=\mu$$
such that either $\lambda^{i+1}\in\supp(\lambda^i)$ with
$\lambda^i$ and $\lambda^{i+1}$ translation
equivalent or
$\lambda^{i+1}\in\supp^2(\lambda^i)$ with
$\lambda^{i+1}=\lambda^i+\epsilon_s-\epsilon_t$ (for some $s$ and $t$)
and $\lambda^i$ and $\lambda^{i+1}$ are
$(\res^{\lambda^i,\pm},\ind^{\lambda^{i+1},\pm})$-translation equivalent.

Suppose that $\lambda,\tilde{\lambda}\in \Lambda_n$ with
$\tilde{\lambda}=\lambda-\epsilon_i+\epsilon_j$, and that $\lambda$
lies in the $(i,j)_-$-hyperplane. By Lemmas \ref{wallcase} and
\ref{samen} we have a bijection $\theta:{\mathcal
  B}(\lambda)\rightarrow{\mathcal B}(\tilde{\lambda})$ which restricts
to a bijection $\theta:{\mathcal B}_n(\lambda)\rightarrow{\mathcal
  B}_n(\tilde{\lambda})$. 
 By Corollary \ref{pmproj}, Theorem
\ref{bigres}, Lemmas \ref{wallcase} and \ref{samen}, and standard properties of
$\ind$ and $\res$ it is clear that weights $\lambda$ and
$\tilde{\lambda}$ are
$(\res^{\lambda,\pm},\ind^{\tilde{\lambda} ,\pm})$ translation
equivalent. It is also easy to see that the adjointness isomorphism is
multiplicative. Thus we can apply Theorem \ref{ripair} and
get a Morita equivalence between the
two blocks $B_n(\lambda)$ and $B_n(\tilde{\lambda})$.

\begin{thm}\label{factran} 
If $\lambda$ and $\mu$ are in the same facet then they are in the same
$(\pm)$-translation class.
\end{thm}
\begin{proof}
By Lemmas \ref{onwall} and \ref{wallcase} it is enough to show that if
$\lambda$ and $\mu$ are in the same facet then there is a chain of
dominant weights
$$\lambda=\lambda^0,\lambda^1,\ldots,\lambda^r=\mu$$ in the same facet
such that $\lambda^{i+1}\in\supp(\lambda^i)$ or
$\lambda^{i+1}\in\supp^2(\lambda^i)$ for each $i$.

Let $x=\lambda+\rho_{\delta}$ and $y=\mu+\rho_{\delta}$, and recall
that in Section \ref{weyl} we associated a function $f$ to each facet
(rather than just a permutation $\pi$ as for an alcove). The proof now
proceeds exactly as for the alcove case (Theorem \ref{alcovemor})
replacing $\pi$ by $f$, until we reach some point where
$f(i)=(k,l)$. In this case we repeatedly add (or subtract) a box from
$y_k$ and subtract (or add) a box to $y_l$ until we reach $x_k$ and
$x_l$. In each of these steps we obtain some
$\lambda^{i+1}=\lambda^i\pm(\epsilon_k-\epsilon_l)\in\supp^2(\lambda)$.
\end{proof}

Applying the results on translation and $(R,I)$-translation
equivalence from Section \ref{towermor} we deduce

\begin{cor}\label{sumup}
If $\lambda$
and $\lambda'$ are in the same facet and $\mu$ and $\mu'$ are such that
$\lambda,\mu\in\Lambda_n$, and $\lambda',\mu'\in\Lambda_m$, and $\mu'$
is the unique weight in ${\mathcal B}(\lambda')$ in the same facet as
$\mu$, then we have:\\ (i)
$$[\Delta_n(\lambda):L_n(\mu)]=[\Delta_m(\lambda'):L_m(\mu')]$$
(ii)
$$\Hom_n(\Delta_n(\lambda),\Delta_n(\mu))\cong
\Hom_m(\Delta_m(\lambda'),\Delta_m(\mu'))$$
(iii) 
$$\Ext^i_n(\Delta_n(\lambda),\Delta_n(\mu))\cong
\Ext^i_m(\Delta_m(\lambda'),\Delta_m(\mu'))$$
for all $i\geq 1$. If further ${\mathcal B}_n(\lambda)$ and ${\mathcal
  B}_m(\lambda')$  contain the same number of simples then the corresponding
blocks are Morita equivalent. 
\end{cor}

By \cite[Theorem 3.4]{dhw} there are always enough local homomorphisms
for the Brauer algebra. Further, by Lemma \ref{onwall} any weight that
is not adjacent to a weight in a less singular facet is translation
equivalent to a weight of smaller total degree. Thus if $\delta>0$
then every weight can be reduced to a weight in some alcove by
translation equivalence and repeated applications of Proposition \ref{alcok}.
This implies

\begin{cor}\label{plus0} If $\delta>0$ then the decomposition numbers
  $[\Delta_n(\lambda):L_n(\mu)]$ for arbitrary $\lambda$ and $\mu$ are
  determined by those for $\lambda$ and $\mu$ in an alcove.
\end{cor}

Note that the restriction on $\delta$ is necessary, as for $\delta<0$
there are no weights in an alcove. In fact for $\delta=-2m$ or
$\delta=-2m+1$ any dominant weight is $\delta$-singular of degree at
least $m$. For the rest of this section we will see what more can be
said in such cases.

We will denote the set of all partitions $\lambda$ with at most $m$
non-zero parts by $\Lambda^{\leq m}$, and the set of those with at
most $m+1$ non-zero parts with $\lambda_{m+1}\leq 1$ by $\Lambda^{\leq
  m,1}$. (Note that $\Lambda^{\leq m}$ is precisely the set of weights
considered in Theorem \ref{fakealc}.) Such weights lie in a union of
facets, but we shall see that together they play a role analogous to
that played by the fundamental alcove in the $\delta>0$ case.

We begin by noting 

\begin{prop}\label{fullrep}
(i) For $\delta=-2m$, every $\delta$-singular weight of degree $m$ is in the
  same block as a unique element of $\Lambda^{\leq m}$.\\ (ii) For
  $\delta=-2m+1$, every $\delta$-singular weight of degree $m$ is in the same
  block as a unique element of $\Lambda^{\leq m,1}$.\\
\end{prop}
\begin{proof}
(i) Suppose
that $\delta=-2m$ and let $\lambda$ be a $\delta$-singular weight of degree
$m$. Then $\lambda+\rho_{\delta}$ is of the form
$$(\ldots,x_1,\ldots,x_2,\ldots,x_m,\ldots,(0),\ldots,-x_m,\ldots,-x_2,
\ldots,-x_1,\ldots,-n,-(n+1),\ldots)$$ where the only elements of
equal modulus are those of the form $\pm x_i$, and the bracketed $0$
may or may not appear. Note that as $\lambda$ is a finite weight the
tail of $\lambda+\rho_{\delta}$ will equal the tail of
$\rho_{\delta}$, i.e. has value $-n$ in position $m+1+n$ for all
$n>>0$, and we assume that this holds for the $n$ in the expression
above (and similar expressions to follow).

First suppose that $\lambda+\rho_{\delta}$ contains $0$. Then $\lambda+\rho_{\delta}$ is
in the same $W_a$-orbit as 
$$\mu+\rho_{\delta}=(x_1,x_2,\ldots,x_m,0,\ldots,-x_m,\ldots,-x_2,\ldots,-x_1,
\ldots,-n,-(n+1),\ldots).$$
Thus the $n-1$ coordinates between the entries $0$ and $-n$ must be
strictly decreasing, which forces
$$\mu+\rho_{\delta}=(x_1,x_2,\ldots,x_m,0,-1,-2,-3,\ldots).$$
Hence we deduce that $\mu\in\Lambda^{\leq m}$ as required.

Next suppose that $\lambda+\rho_{\delta}$ does not contain $0$. Then
there are two cases depending on the parity of the number of positive
entries in $\lambda+\rho_{\delta}$. The first case is when
$\lambda+\rho_{\delta}$ is in the same $W_a$-orbit as
$$\mu+\rho_{\delta}=(x_1,x_2,\ldots,y,\ldots,x_m,\ldots,-x_m,\ldots,-x_2,\ldots,-x_1,
\ldots,-n,-(n+1),\ldots)$$ where $y$ is some positive integer and all
entries after $x_m$ are negative.  Arguing as above we see that
$$\mu+\rho_{\delta}=(x_1,x_2,\ldots,y,\ldots,x_m,-1,-2,-3,\ldots).$$
But this vector is $\delta$-singular of degree $m+1$, which contradicts our
assumptions on $\lambda$.

The second case is when $\lambda+\rho_{\delta}$ is in the same
$W_a$-orbit as 
$$\mu+\rho_{\delta}=(x_1,x_2,\ldots,x_m,\ldots,-x_m,\ldots,-x_2,\ldots,-x_1,
\ldots,-n,-(n+1),\ldots)$$ 
where all entries after $x_m$ are negative.
But this implies that $\mu+\rho_{\delta}$ has $n$ strictly decreasing
coordinates between the entries $0$ and $-n$, which is impossible.

The argument for (ii) is very similar. We see that $\lambda+\rho_{\delta}$ is
in the same $W_a$-orbit as either
$$\mu+\rho_{\delta}=(x_1,x_2,\ldots,x_m,\ldots,-x_m,\ldots,-x_2,\ldots,-x_1,
\ldots,-n-\frac{1}{2},-(n+1)-\frac{1}{2},\ldots)$$
or
$$\mu+\rho_{\delta}=(x_1,x_2,\ldots,y,\ldots,x_m,\ldots,-x_m,\ldots,-x_2,\ldots,-x_1,
\ldots,-n-\frac{1}{2},-(n+1)-\frac{1}{2},\ldots)$$
where in each case all entries after $x_m$  are negative.

In the first case we deduce as above that 
$$\mu+\rho_{\delta}=(x_1,x_2,\ldots,x_m,-\frac{1}{2},-\frac{3}{2},-\frac{5}{2},\ldots)$$
and so $\mu\in\Lambda^{\leq m}\subset \Lambda^{\leq m,1}$ as
required. In the second case, as $\mu+\rho_{\delta}$ must be
$\delta$-singular of degree $m$, we deduce that
$$\mu+\rho_{\delta}=(x_1,x_2,\ldots,y,\ldots, x_m,-\frac{1}{2},-\frac{3}{2},-\ldots,
-y+1,\hat{y},-y-1\ldots)$$
where $\hat{y}$ denotes the omission of the entry $y$. But this
element is in the same $W_a$-orbit as 
$$\nu+\rho_{\delta}=(x_1,x_2,\ldots,x_m,\frac{1}{2},-\frac{3}{2},-\frac{5}{2},
-\frac{7}{2},\ldots)$$
(by swapping $y$ and $-\frac{1}{2}$ with a change of signs, and rearranging
to get a decreasing sequence). Thus $\lambda$ is in the same
$W_a$-orbit as $\nu$, and $\nu=(\nu_1,\ldots,\nu_m,1)
\in\Lambda^{\leq m,1}$ as required.
\end{proof}

Although the weights in $\Lambda^{\leq m}$ (respectively in
$\Lambda^{\leq m,1}$) lie in several different facets, the next result
shows that all these facets have equivalent representation
theories.

\begin{prop}\label{bigalc}
Let $\delta=-2m$ (respectively $\delta=-2m+1$) and $\lambda\in
\Lambda^{\leq m}$ (respectively $\lambda\in\Lambda^{\leq m,1}$). Then
$\lambda$ and $\lambda'$ are translation equivalent if and only if 
$\lambda'\in\Lambda^{\leq m}$ (respectively $\lambda'\in\Lambda^{\leq m,1}$).
\end{prop}
\begin{proof}
Note that all weights in $\Lambda^{\leq m}$ (respectively in
$\Lambda^{\leq m,1}$) are $\delta$-singular of degree $m$, and that
any pair of such weights can be linked by a chain of weights in the
same set differing at each stage only by the addition or subtraction
of a single block. By Lemma \ref{onwall} we see that any such pair is
translation equivalent.

For the reverse implication, we first consider the $\delta=-2m$ case,
with $\lambda\in \Lambda^{\leq m}$, and suppose that
$\lambda'\in\supp(\lambda)$ is not an element of $\Lambda^{\leq m}$.
Then we must have
$\lambda'=(\lambda_1,\lambda_2,\ldots,\lambda_m,1)$. Now
$x=\lambda+\rho_{\delta}$ and $x'=\lambda'+\rho_{\delta}$ differ only
in the $m+1$st coordinate, which is $0$ or $1$ respectively. If $f$ is
the function associated to the facet containing $x$ and $f'$ is the
corresponding function for $x'$ then the only difference is that
$f(1)=m+1$ while $f'(1)=(m+1,m+2)$. Thus
$\lambda'\in\overline{F}\backslash F$, where $F$ is the facet
containing $\lambda$, and so by Lemma \ref{onwall} this pair cannot be
translation equivalent.

The case $\delta=-2m+1$ is similar. Arguing as above we have that
$\lambda'=(\lambda_1,\ldots,\lambda_m,1,1)$ or
$\lambda'=(\lambda_1,\ldots,\lambda_m,2)$, and in each case it is easy
to show that the pair $\lambda$ and $\lambda'$ are not translation
equivalent.
\end{proof}

Combining Propositions \ref{fullrep} and \ref{bigalc} we deduce that
for $\delta<0$ and $\lambda,\lambda'\in\Lambda^{\leq m}$ (respectively
$\lambda,\lambda'\in\Lambda^{\leq m,1}$) there is a bijection $\theta:{\mathcal
  B}(\lambda)\rightarrow {\mathcal B}(\lambda')$ which as before we
will denote by $\theta(\mu)=\mu'$.  Applying the results from Section
\ref{towermor} we obtain

\begin{cor} Let $\delta<0$ and  $\lambda, \lambda'\in\Lambda^{\leq m}$
  (respectively $\lambda,\lambda'\in\Lambda^{\leq m,1}$). If
  $\lambda,\lambda'\in\Lambda_n$ and $\lambda',\mu'\in\Lambda_l$ then
  (i--iii) of Corollary \ref{sumup} hold. If further ${\mathcal
    B}_n(\lambda)$ and ${\mathcal B}_l(\lambda')$ contain the same
  number of elements then the corresponding blocks are Morita equivalent.
\end{cor}

As in Corollary \ref{plus0}, we obtain the following application of
Proposition \ref{alcok}.

\begin{cor}\label{neg0} If $\delta<0$ then the decomposition numbers
  $[\Delta_n(\lambda):L_n(\mu)]$ for arbitrary $\lambda$ and $\mu$ are
  determined by those for $\lambda$ and $\mu$ in a singular facet of
  degree $m$.
\end{cor}

Combining Corollaries \ref{plus0} and \ref{neg0} with our earlier
remarks we obtain

\begin{thm} \label{need0} For $\delta\in\ZZ$ non-zero
the decomposition numbers $[\Delta_n(\lambda):L_n(\mu)]$ for arbitrary
$\lambda$ and $\mu$ are determined by those for $\lambda$ and $\mu$ in
${\mathcal B}(0)$.
\end{thm}

Thus (at least at the level of decomposition numbers) is it enough to
restrict attention to a single block of the Brauer algebra.

\begin{rem}
The decomposition numbers for the module $\Delta_n(0)$ are known by
\cite[Proposition 5.1 and Theorem 5.2]{cdm}.
\end{rem}

We would also like the representation theory to be independent of
$\delta\in\ZZ$, in the sense that it should depend only on the geometry of
facets. For weights in alcoves, this would in large part follow if we
could show that decomposition numbers are given by some kind of
parabolic Kazhdan-Lusztig polynomials. In the remaining sections we will
consider some evidence for this.

\section{Block graphs for the Brauer algebra}\label{manygraphs}

Recall from Section \ref{basics} the definition of a maximal
balanced partition.  Let $\MBS_{\delta}(\lambda)$ be the
directed graph with vertex set $V_{\delta}(\lambda)$ and edge
$\mu\rightarrow\tau$ if $\mu$ is a maximal $\delta$-balanced
subpartition of $\tau$.

The above graph appears to depend both on $\lambda$ and $\delta$,
while the alcove geometry associated to $W_a$ does not. 
Let $\Alc$ be the directed graph with vertex set the set of alcoves
for $W_a$ in $A^+$, and an edge $A\rightarrow B$ if the
closures of $A$ and $B$ meet in a hyperplane and this hyperplane
separates $A_0$ and $B$. (Note that the former condition corresponds
to $B=(ij)_-A$ for some reflection $(ij)_-$.)

Our goal in this section is to show that all the graphs
$\MBS_{\delta}(\lambda)$ are in fact isomorphic, and are isomorphic to
the alcove graph $\Alc$. This will be our first evidence that 
the representation theory depends only on the geometry of
facets.

Recall that for $\lambda \in X^+$ we have $\lambda + \rho_\delta\in
A^+$ (the set of strictly decreasing sequences in $\mathbb{E}^\infty$)
and the $\delta$-dot action of $W_a$ on $\lambda$ corresponds
to the usual action of $W_a$ on $\lambda + \rho_\delta$. For
the rest of this section we will work with the usual action of
$W_a$ on $A^+$.

For $v\in A^+$ we define 
$$V(v)=W_av\cap A^+.$$ We define a partial order on $A^+$ by
setting $x\leq y$ if $y-x\in \mathbb{E}^f$ and all entries in $y-x$
are non-negative. For $v\in A^+$ we define a directed graph $\Graph(v)$
with vertex set $V(v)$ and arrows given as follows. If $x,y\in V(v)$,
we set $x\rightarrow y$ if and only if $x<y$ and there is no $z\in
V(v)$ with $x<z<y$. The reason for introducing this graph is clear
from

\begin{prop}\label{iso1}
For $\lambda \in X^+$ we have $\MBS_\delta (\lambda^T) \cong \Graph(\lambda
+\rho_\delta)$.
\end{prop}
\begin{proof}
By Theorem \ref{geoblock} we have a bijection between
$V_\delta(\lambda^T)$ and $V(\lambda + \rho_\delta)$. Moreover, for
$\mu^T, \nu^T, \tau^T \in V_\delta (\lambda^T)$ we have $\mu^T \subset
\nu^T \subset \tau^T$ if and only if $\mu + \rho_\delta < \nu +
\rho_\delta < \tau +\rho_\delta$. Thus the two graph structures on
these vertex sets are preserved under the correspondence.
\end{proof}

Recall the definition of singletons from Section \ref{weyl}. Define
$v_{reg}$ to be the subsequence of $v$ consisting only if its
singletons. For example, if $v$ begins
$(9,8,7,0,-1,-2,-7,-9,-11,\ldots)$ then $v_{reg}$ begins
$(8,0,-1,-2,-11,\ldots).$ Note that if $v\in A^+$ then $v_{reg}\in
A^+$ and $|(v_{reg})_i|\neq |(v_{reg})_j|$ for all $i\neq j$.
Therefore $v_{reg}$ is a regular element in $\mathbb{E}^\infty$ (as it
does not lie on any reflecting hyperplane). We define the
\emph{regularisation map} $\Reg \, :\, A^+ \longrightarrow A^+$ by
setting $$\Reg(v) = v_{reg}.$$ The key result about the regularisation
map is

\begin{prop}\label{iso2}
For all $v\in A^+$ we have $$\Graph(v) \cong \Graph(v_{reg}).$$
\end{prop}
\begin{proof}
We first observe that the map $\Reg$ gives rise to a bijection between
$V(v)$ and $V(v_{reg})$. For the set of doubletons is an invariant of the
elements in $V(v)$, and given this set there is a unique way of
adding the doubletons into an element of $V(v_{reg})$ keeping the
sequence strictly decreasing.  Now suppose that $x,y\in A^+$ and $a\in
\mathbb{R}$ are such that
$$s=(x_1, ... , x_i, a, x_{i+1}, ...)\in A^+ \quad
\wand\quad t=(y_1, ... , y_j, a, y_{j+1}, ...)\in A^+.
$$ Then it is easy to see that $x<y$ if and only if $s<t$.  However, this
implies that the set of edges coincide under the map $\Reg$, as
required.
\end{proof}

\begin{cor} For all $v, v'\in A^+$ we have 
$$\Graph(v) \cong\Alc.$$
Hence for all $\delta,\delta'\in\ZZ$ and
  $\lambda,\lambda'\in X^+$ we have
$$\MBS_{\delta}(\lambda)\cong \MBS_{\delta'}(\lambda').$$
\end{cor}

\begin{proof} 
Note that any $v\in \Reg(A^+)$ lies inside an alcove. For any vector
$v\in\Reg(A^+)$ the maximal weights below $v$ in the same orbit lie in
the alcoves below and adjacent to the alcove containing $v$.
Thus for $v\in
\Reg(A^+)$ it is clear that we have $\Graph(v)\cong \Alc$. Now
the result follows for general $v$ from  Proposition \ref{iso1}, and
in its $\MBS$ form from
Proposition \ref{iso2}.
\end{proof}

It will be convenient to give $\Alc$ the structure of a graph with
coloured edges. An edge in $\Alc$ corresponds to reflection from an
alcove $A$ to an alcove $B$ through the facet separating them. The
action of $W_a$ on weights induces a corresponding action on facets,
and we shall say that two edges have the same colour if and only if
the corresponding facets lie in the same orbit.

We conclude this section with one final graph $\Par^+_e$ isomorphic to
$\Alc$, whose structure can be described explicitly.

We fix the element $v=(-1,-2,-3,-4,...)\in A^+$. Using the action of
$W_a$ we can see that every $x\in V(v)$ corresponds uniquely
to a strictly decreasing partition with an even number of parts,
obtained by ignoring all parts of $x$ which are negative. For example,
the element $x=(6,5,3,1,-2,-4,-7,-8,...)$ corresponds to $(6,5,3,1)$
while $v$ corresponds to $\emptyset$. Thus if we write $P^+_e$ for
the set of strictly decreasing partitions with an even number of parts
then we have a bijection
$$\phi \, : \, V(v) \longrightarrow P^+_e.$$

Consider the usual partial order $\subseteq$ on $P^+_e$ given by
inclusion of partitions (viewed as Young diagrams). It is clear that
the partial order $\leq$ on $V(v)$ corresponds to the partial order
$\subseteq$ on $P^+_e$ under the bijection $\phi$. Define a graph
$\Par^+_e$ with vertex set $P^+_e$ and an arrow $\lambda \rightarrow
\mu$ if and only if $\lambda \subset \mu$ and there is no $\nu\in
P^+_e$ with $\lambda \subset \nu \subset \mu$. It is easy to verify
that the map $\phi$ induces
a graph isomorphism between $\Graph(v)$ and $\Par^+_e$.

The graph $\Par^+_e$ can easily be described explicitly as follows. For
$\lambda, \mu\in P_e^+$, there is an arrow $\lambda \rightarrow \mu$
if and only if either
\begin{equation}\label{cond1}
\lambda = (\lambda_1, \ldots , \lambda_n)\,\,\quad\wand\quad
\mu=(\lambda_1, \ldots , \lambda_{i-1}, \lambda_i + 1,
\lambda_{i+1}, \ldots, \lambda_n)
\end{equation}
 or
\begin{equation}\label{cond2}
\quad \lambda =(\lambda_1, \ldots , \lambda_n) \ \mbox{\rm
  with}\
\lambda_n\geq 3 \quad\wand\quad
\mu= (\lambda_1,
\ldots , \lambda_n , 2,1).
\end{equation}
 To see this, first observe that in both cases there is no $\nu\in
 P_e^+$ with $\lambda \subset \nu \subset \mu$. Moreover, if $\lambda,
 \mu\in P^+_e$ with $\lambda \subset \mu$ then $\mu$ can be obtained
 from $\lambda$ by applying (\ref{cond1}) and (\ref{cond2}) repeatedly.

\section{Coxeter systems and parabolic Kazhdan-Lusztig polynomials}
\label{klp}

In this section we will introduce parabolic Kazhdan-Lusztig
polynomials associated with the pair $(W_a,W)$. We briefly
review the relevant theory; details can be found in \cite{humcox} and
\cite{soergel1}.

\begin{figure}[ht]
\includegraphics{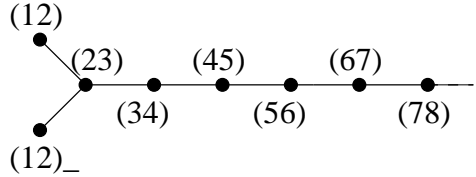}
\caption{The type $D_{\infty}$ Coxeter system}
\label{DCox}
\end{figure}

Recall that a Coxeter system is a pair $(G,S)$ consisting of a group
$G$ and a set $S$ of generators of $G$ such that all relations in $G$
are of the form 
$$(ss')^{m(s,s')}=1$$ where $m(s,s)=1$ and $m(s,s')=m(s',s)\geq 2$
otherwise (including the possibility that $m(s,s')=\infty$ denoting no
relation between $s$ and $s'$). Note that the group $G$ does not need
to be finite (although this is often assumed). Given a Coxeter system,
the associated Coxeter graph is the graph with vertices the elements of
$S$, and $m(s,s')-2$ edges between $s$ and $s'$ (or no edges when
$m(s,s')=\infty$). For example the $D_{\infty}$ Coxeter graph is given
by the graph shown in Figure \ref{DCox}.

It is easy to verify that our group $W_a$ is generated by the
elements $\{(12)_-,(i\,i+1): i\geq 1\}$ and satisfies the relations
given by the Coxeter graph in Figure \ref{DCox}. Thus it must be a
quotient of the Coxeter group of type $D_{\infty}$. However, for each
choice of $n$, the subsystem generated by the first $n$ generators is
precisely the type $D_{n}$ system (see \cite[Planche IV]{bourbaki}),
and so there can be no further relations, and our group is the type
$D_{\infty}$ Coxeter group associated to the given generators.

Given a Coxeter system $(G,S)$, any subgroup $G'$ generated by a
subset $S'$ of $S$ defines a parabolic subsystem $(G',S')$. In our
case the group $W$ clearly arises in this way from the generators of
the form $(i\,i+1)$ and so is a type $A_{\infty}$ parabolic subgroup of
$W_a$. 

When $\delta\geq 0$ there is a bijection from $W_a$ to the set
of alcoves, given by $w\longmapsto w.0$. We will henceforth identify
elements of $W_a$ with alcoves via this map.  Under this
bijection the standard length function on our Coxeter system
associated to $W_a$ (given in terms of the number of terms
occurring in a reduced expression for $w$) corresponds to the number of
reflection hyperplanes between $0$ and $w\cdot_{\delta}0$.

We define $W^{a}$ to be the subset of $W_a$ corresponding
to the alcoves in $X^+$. By Proposition \ref{Wis}(ii) we then have a
bijection
$$W\times W^{a}\rightarrow W_a.$$ We are thus in a position to define
$D_{\infty}/A_{\infty}$ parabolic Kazhdan-Lusztig polynomials
following \cite{deodhar} (although we use the notation of
\cite[Section 3]{soergel1}). Their precise definition and general
properties need not concern us, instead we will give a recursive
construction corresponding to stepping away from the root of
$\Alc$. To do this we will first need to define a partial order on
weights (and alcoves).

Two weights $\lambda$ and $\mu$ such that $\mu=w.\lambda$ for some
reflection $w$, lie in different components of the space formed by
removing this hyperplane.  We say that $\lambda<\mu$ if $\lambda$ is
in the component containing the fundamental alcove. This extends to
give a partial order on weights, which in turn induces a partial order
on alcoves. This agrees with the path-from-root order on $\Alc$. Two
alcoves are said to be \emph{adjacent} if there is precisely one
reflecting hyperplane between them (i.e. they are adjacent in $\Alc$).

Suppose that $\nu$ and $\mu$ are dominant weights in adjacent alcoves
with $\nu=s\cdot_{\delta}\mu>\mu$. Given a dominant weight $\lambda \in
W_a\cdot_{\delta}\mu$ we define $\kappa_{\lambda}(\nu,\mu)$ to
be the unique weight such that 
$$(\kappa_{\lambda}(\nu,\mu),\lambda)=
(w\cdot_{\delta}\mu,w\cdot_{\delta}\nu)$$
i.e. $(\kappa_{\lambda}(\nu,\mu),\lambda)$ is an edge of the same
colour as $(\mu,\nu)$ in $\Alc$.

We next define certain polynomials $n_{\nu,\lambda}$ (in an
indeterminate $v$) for regular
weights $\lambda$ and $\mu$ in the following recursive manner. Let
$e_{\lambda}$ as $\lambda$ runs over the  regular weights be a set of
formal symbols.

\begin{enumerate}
\item[(i)]  We set
$n_{\nu,\lambda}=0$ if $\lambda\not\leq\nu$ or $\lambda\notin
W_a\cdot_{\delta}\nu$ or either $\lambda$ or $\nu$ is
non-dominant.
\item[(ii)] We set $n_{0,0}=1$ and $N(0)=e_{0}$.
\item[(iii)] For each $\nu>0$ regular dominant, there exists some $\mu$ regular
  dominant below it such that $\mu=s\cdot_{\delta}\nu$ and $(\nu,\mu)$
  are in adjacent alcoves. Pick any such $\mu$. Then for any dominant
  $\lambda$ with $\lambda=w\cdot_{\delta}\nu$ for some $w$ and
  $\kappa=\kappa_{\lambda}(\nu.\mu)$ we set
$$\hat{n}_{\nu,\lambda}=\pr_+(\kappa)\left(n_{\mu,\kappa}
  +v^{l(\kappa)-l(\lambda)}n_{\mu,\lambda}\right)$$ where
  $\pr_+(\kappa)=1$ if $\kappa\in X^+$ and $\pr_+(\kappa)=0$
  otherwise. Note that for $\kappa\in X^+$ we have
  $l(\kappa)-l(\lambda)=-1$ if $\kappa<\lambda$, respectively $+1$ if
  $\kappa>\lambda$.  Let $\hat{N}(\nu)$ be the sum
$$\hat{N}(\nu)=\sum_{\lambda}\hat{n}_{\nu,\lambda}e_{\lambda}.$$
and $R(\nu)$ be the set of $\lambda<\nu$ such that
$\hat{n}_{\nu,\lambda}(0)\neq 0$. Then 
$$N(\nu)=\hat{N}(\nu)-\sum_{\lambda\in R(\nu)}\hat{n}_{\nu,\lambda}(0)N(\lambda)$$
and $n_{\nu,\lambda}$ is the coefficient of $e_{\lambda}$ in $N(\nu)$.
\end{enumerate}

It is a consequence of (Deodhar's generalisation of) Kazhdan-Lusztig
theory that this process is well defined (so does not depend on the
choice of $\mu$ in step (iii)), and that each $n_{\nu,\lambda}$ is a
polynomial in $v$ with $n_{\nu,\lambda}(0)\neq 0$ only if $\lambda=\nu$.

\section{Some low rank calculations for $\delta=1$}
\label{klpoly}

To illustrate the various constructions so far, we will consider the
case when $\delta=1$, and examine the regular block containing the
weight $0$. First we calculate the associated parabolic
Kazhdan-Lusztig polynomials, and then we compare these with the
representation theoretic results.

We will also need to consider the block containing $(1)$. As this is
in the same alcove as $0$ these two blocks are translation
equivalence. However, in this simple case case we do not obtain any
simplification to the calculations by applying the results from
Section 6; instead the results can be considered as a verification of
the general theory in this special case.

\begin{figure}[ht]
\includegraphics{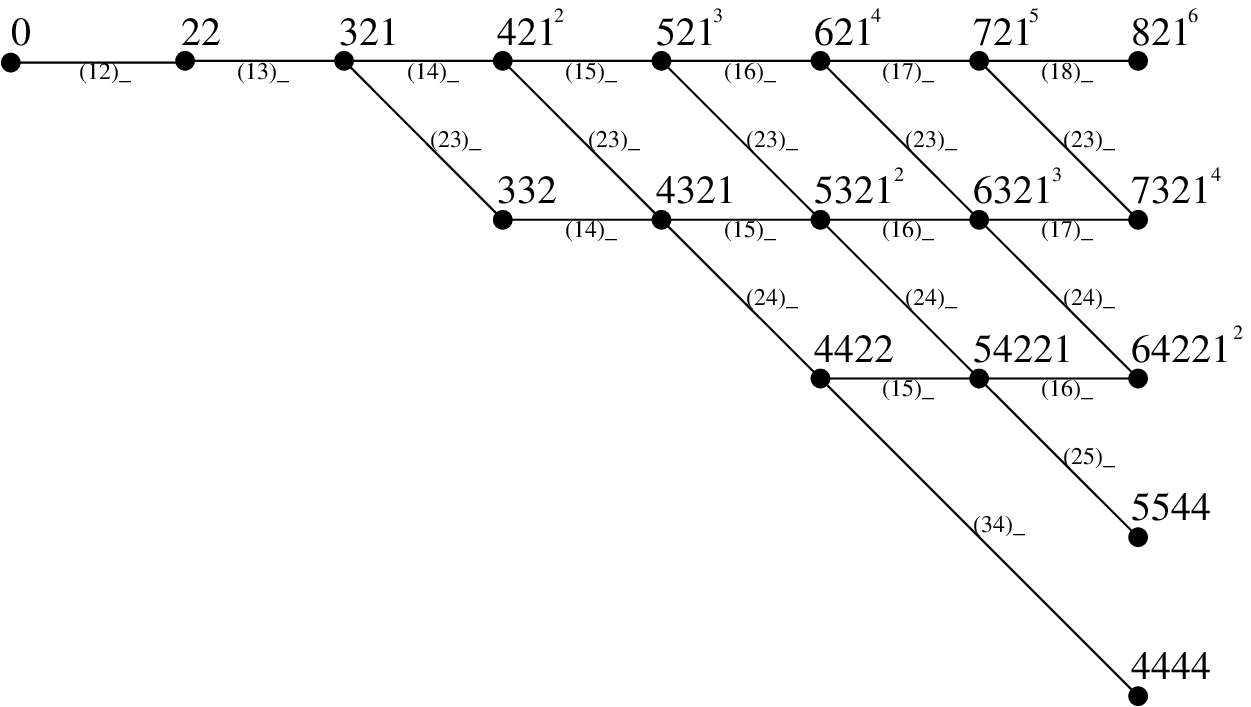}
\caption{The block of the weight $0$ for $\delta=1$, up to degree 16}
\label{blockex}
\end{figure}

In Figure \ref{blockex} we have listed all dominant weights of degree
at most 16 that are in the same block as the weight $0$.  We will
abbreviate weights in the same manner as partitions (and so write for
example $(1^3)$ instead of $(1,1,1)$). An edge between two weights indicates
that they are in adjacent alcoves, and the label $(ij)_-$ corresponds to
the reflection hyperplane between them. (Clearly only weights of the
form $(ij)_-$ can arise as such labels.)

Given this data we can now compute the $n_{\lambda,\mu}$. The final
results are shown in Figure \ref{kltable}. We start with the weight
$0$ having $n_{(0),(0)}=1$. Reflecting through $(12)_-$ we obtain the
weight $(22)$, and we see that $n_{(22),(0)}=v$. (Note that the term
$n_{\lambda,\lambda}$ is always $1$.) Continuing we reflect $(22)$
through $(13)_-$ to obtain $(321)$. As $(12)_-(13)_-\cdot_{\delta}0
=(2,-1,3)$ is not dominant we see that the only non-zero term apart
from $n_{(321),(321)}$ is $n_{(321),(22)}=v$. Identical arguments give
all polynomials $n_{\nu,\lambda}$ where $\nu$ is on the top row of
Figure \ref{blockex}. For the second row we obtain four terms as $\nu$
and $\mu$ both give dominant weights under the action of
$(1i)_-(23)_-$ for suitable $i$ (as the parallelogram with $\nu$ as
highest term has identically labelled parallel sides).

For $(4422)$ we must observe that
$$(14)_-(23)_-(24)_-\cdot_{\delta}(4321)=(22)$$ and similar results
give the remaining cases. In all of these cases we have no constant
terms arising at any stage (apart from in $n_{\lambda,\lambda}$), and
hence $N(\nu)=\hat{N}(\nu)$ for every weight considered.

\begin{figure}\label{klcalc}
\begin{tabular}{l|cccccccccccccccccc}
$\mu\quad\ \ \lambda$ & \begin{sideways}$0$\end{sideways} 
& \begin{sideways}$22$\end{sideways}
    & \begin{sideways}$321$\end{sideways}  & 
\begin{sideways}$4211$\end{sideways}  & 
\begin{sideways}$332$\end{sideways}
  & \begin{sideways}$521^3$\end{sideways}  & 
\begin{sideways}$4321$\end{sideways} &  
\begin{sideways}$621^4$\end{sideways}
  &\begin{sideways}$53211$\end{sideways}  &
\begin{sideways}$4422$\end{sideways} &
\begin{sideways}$721^5$\end{sideways}
  & \begin{sideways}$6321^3$\end{sideways}  &
\begin{sideways}$54221$\end{sideways} 
      &
\begin{sideways}$821^6$\end{sideways}  &
\begin{sideways}$7321^4$\end{sideways} &
\begin{sideways}$64221^2$\end{sideways} &
\begin{sideways}$552^3$\end{sideways} &
\begin{sideways}$4^4$\end{sideways} \\
\hline
$ 0$ &$ 1$&&&&&&&&&&&&&&&&&\\
$22$ & $v$&$1$&&&&&&&&&&&&&&&&\\
$321$ & .&$v$&$1$&&&&&&&&&&&&&&&\\
$4211$ & .&.&$v$&$1$&&&&&&&&&&&&&&\\
$332$ & .&.&$v$&.&$1$&&&&&&&&&&&&&\\
$521^3$ & .&.&.&$v$&.&$1$&&&&&&&&&&&&\\
$4321$& .&.&$v^2$&$v$&$v$&.&$1$&&&&&&&&&&&\\
$621^4$ &.&.&.&.&.&$v$&.&$1$&&&&&&&&&&\\
$53211$ &.&.&.&$v^2$&.&$v$&$v$&.&$1$&&&&&&&&&\\
$4422$ &.&$v^2$&$v$&.&.&.&$v$&.&.&$1$&&&&&&&&\\
$721^5$ & .&.&.&.&.&.&.&$v$&.&.&$1$&&&&&&&\\
$6321^3$ &.&.&.&.&.&$v^2$&.&$v$&$v$&.&.&$1$&&&&&&\\
$54221$ &.&.&.&.&.&.&$v^2$&.&$v$&$v$&.&.&$1$&&&&&\\
$821^6$ &.&.&.&.&.&.&.&.&.&.&$v$&.&.&$1$&&&&\\
$7321^4$&.&.&.&.&.&.&.&$v^2$&.&.&$v$&$v$&.&.&$1$&&&\\
$64221^2$&.&.&.&.&.&.&.&.&$v^2$&.&.&$v$&$v$&.&.&$1$&&\\
$552^3$&.&.&.&.&$v^2$&.&$v$&.&.&.&.&.&$v$&.&.&.&$1$&\\
$4^4$&$v^2$&$v$&.&.&.&.&.&.&.&$v$&.&.&.&.&.&.&.&$1$\\
\end{tabular}
\caption{The polynomials $n_{\lambda,\mu}$ for $\delta=1$ and
  $|\lambda|\leq 16$}
\label{kltable}
\end{figure}

Next we will determine the structure of certain low rank standard
modules for $B_n(1)$ in the block containing $0$. These will then be
compared with the Kazhdan-Lusztig polynomials calculated above.  We
will proceed in stages, and will also need to consider the structure
of modules in the block containing $(1)$. The submodule structure of
modules will be illustrated diagrammatically, where a simple module
$X$ is connected by a line to a simple module $Y$ above it if there is
a non-split extension of $X$ by $Y$. Note that in this section we
follow the usual labelling of modules by partitions (as in \cite{cdm})
and \emph{not} via the transpose map by weights.

\subsection{The case $n\leq 6$}

When $n=0$ or $n=2$ we have 
$$\Delta_n(0)=L_n(0)$$
by quasihereditary (and the absence of any other simples in the same block).

When $n=4$ we have
$$\Delta_4(22)=L_4(22)\quad\quad\wand\quad\quad \Delta_4(0)=
\xymatrix@=10pt{ L_4(0) \ar@{-}[d]\\
L_4(22)}$$
by quasi-hereditary and (\ref{twoboxstuff}).

When $n=6$ we have 
$$\Delta_6(321)=L_6(321)\quad\quad\wand\quad\quad \Delta_6(22)=
\xymatrix@=10pt{ L_6(22) \ar@{-}[d]\\
L_6(321)}
$$ as in the case $n=4$. For the remaining module $\Delta_6(0)$ we
know that $[\Delta_6(0):L_6(22)]=1$ by localising to $n=4$. Applying
Proposition \ref{projelim} with $\mu=(32)$ (as this weight is minimal
in its block) we see that $L_6(321)$ cannot occur in $\Delta_6(0)$.
Hence we have that
$$ \Delta_6(0)=
\xymatrix@=10pt{ L_6(0) \ar@{-}[d]\\
L_6(22).}
$$
The odd $n$ cases are very similar. Arguing as above we see that
$$\Delta_1(1)=L_1(1),\quad\quad \Delta_3(21)=L_3(21),\quad\quad
\Delta_3(1)=
\xymatrix@=10pt{ L_3(1) \ar@{-}[d]\\
L_6(21)}
$$
and for $n=5$ that
$$\Delta_5(311)=L_5(311),\quad\quad \Delta_5(21)=
\xymatrix@=10pt{ L_5(21) \ar@{-}[d]\\
L_5(311)},\quad\quad
\Delta_5(1)=
\xymatrix@=10pt{ L_5(1) \ar@{-}[d]\\
L_5(21).}$$

\subsection{The case $n=7$}

As above, we deduce from quasi-hereditary and
(\ref{twoboxstuff}) that
$$\Delta_7(41^3)=L_7(41^3),\quad\quad
\wand
\quad\quad \Delta_7(311)=
\xymatrix@=10pt{ L_7(311) \ar@{-}[d]\\
L_5(41^3).}$$

For the remaining two standard modules, all composition multiplicities
are known (by localising to the case $n=5$) except for those for the
\lq new' simple $L_7(41^3)$. However, this does not
occur in $\Delta_7(21)$ or in $\Delta_7(1)$ by an application
of Proposition \ref{projelim} with $\mu=(41^2)$ (as this weight is
minimal in its block). Thus we have that 
$$ \Delta_7(21)=
\xymatrix@=10pt{ L_7(21) \ar@{-}[d]\\
L_8(311)}\quad\quad\wand\quad\quad
 \Delta_7(1)=
\xymatrix@=10pt{ L_7(1) \ar@{-}[d]\\
L_7(21).}$$

\subsection{The case $n=8$}

This is similar to the preceding case. 
We have that
$$\Delta_8(421^2)=L_8(421^2),\quad\quad
\Delta_8(332)=L_8(332),\quad\quad\wand
\quad\quad \Delta_8(321)=\!\!\!\!
\xymatrix@!0@R=30pt@C=30pt{ &L_8(321) \ar@{-}[dl]\ar@{-}[dr]&\\
L_8(421^2)&& L_8(332).}
$$

For the remaining two standard modules, all composition multiplicities
are known (by localising to the case $n=6$) except for those involving
$L_8(421^2)$ and $L_8(332)$. However, neither of
these occurs in $\Delta_8(22)$ or in $\Delta_8(0)$ by an application
of Proposition \ref{projelim} with $\mu=(322)$, respectively
$\mu=(3211)$. Thus we have that 
$$ \Delta_8(22)=
\xymatrix@=10pt{ L_8(22) \ar@{-}[d]\\
L_8(321)}\quad\quad\wand\quad\quad
 \Delta_8(0)=
\xymatrix@=10pt{ L_8(0) \ar@{-}[d]\\
L_8(22).}$$

\subsection{The case $n=9$}
In this case we have $6$ standard modules, labelled by $(1)$, $(21)$,
$(311)$, $(41^3)$, $(51^4)$ and $(3^3)$. By (\ref{twoboxstuff}) there
is a homomorphism from each of these standards to the preceding one in
the list, except in the case of $(3^3)$. For this weight we instead
use (\ref{squarestuff}) which tells us that
$$[\Delta_9(311):L_9(3^3)]=1.$$

As in earlier cases, we have that
$$\Delta_9(51^4)=L_9(51^4),\quad\quad
\Delta_9(3^3)=L_9(3^3),\quad\quad\wand
\quad\quad \Delta_9(41^3)=
\xymatrix@=10pt{ L_9(41^3) \ar@{-}[d]\\
L_9(3^3).}$$
The module $L_9(51^4)$ cannot occur in any other standards, by applying
Proposition \ref{projelim} with $\mu=(41^4)$, and similarly $L_9(3^3)$
can only occur in $\Delta_9(311)$, by taking $\mu=(331)$. By the
above observations and localisation to $n=7$ we deduce that
$$
\Delta_9(311)=\!\!\!
\xymatrix@!0@R=30pt@C=30pt{ &L_9(311) \ar@{-}[dl]\ar@{-}[dr]&\\
L_9(41^3)&& L_9(3^3)}\quad\quad
\Delta_9(21)=
\xymatrix@=10pt{ L_9(21) \ar@{-}[d]\\
L_9(311)}\quad\quad
\Delta_9(1)=
\xymatrix@=10pt{ L_9(1) \ar@{-}[d]\\
L_9(21).}
$$

We will need to consider $\res_{10}\Delta_{10}(321)$. For this we need
to understand the various standard modules arising in the short exact
sequence (\ref{resrule}) in this case. First note that $(32)$ and
$(221)$ are the unique weights in their respective blocks when
$n=9$. For the weights $(331)$, $(322)$ and $(3211)$ there is exactly
one larger weight in the same block in each case, respectively
$(4311)$, $(4221)$, and $(3321)$.

It follows from the above remarks, (\ref{twoboxstuff}), and (\ref{resrule}) 
that $\res_{10}\Delta_{10}(321)$ has a short
exact sequence
\begin{equation}
\label{9res}
0\too A\too \res_{10}\Delta_{10}(321)\too B\too 0
\end{equation}
where 
\begin{equation}
\label{9lower}
A\cong L_9(221)\oplus \!\!
\xymatrix@!0@R=15pt@C=20pt{ &L_9(311)&\\
L_9(41^3)&& L_9(3^3)}\!\!\oplus
L_9(32)
\end{equation}
and
\begin{equation}
\label{9upper}
B\cong\! \xymatrix@!0@R=15pt@C=20pt{L_9(421)\\
L_9(432)}\!\oplus\!
\xymatrix@!0@R=15pt@C=20pt{L_9(331)\\
L_9(4311)}\!\oplus\!
\xymatrix@!0@R=15pt@C=20pt{L_9(322)\\
L_9(4221)}\!\oplus\!
\xymatrix@!0@R=15pt@C=20pt{L_9(3211)\\
L_9(3321).}
\end{equation}

\subsection{The case $n=10$}
From now on, we will summarise the results obtained for each value of
$n$ in a single diagram, together with an explanation of how they were
derived. In each such diagram we shall illustrate the structure of
individual modules as above, but label simple factors just by the
corresponding partition. We will indicate the existence of a
homomorphism between two modules by an arrow. (It will be clear which
standard module is which by the label of the simple in the head.) 

\begin{figure}
\includegraphics{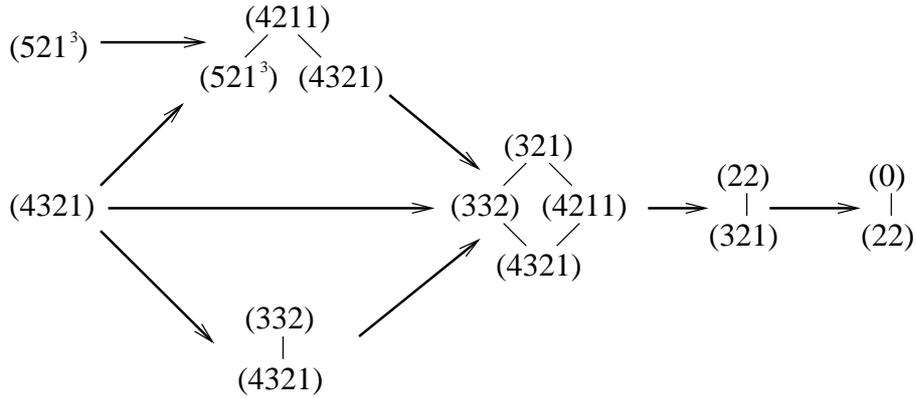}
\caption{The block containing $0$ when $n=10$}\label{case10}
\end{figure}

For $n=10$ we claim that the structure of the block containing $(0)$
is given by the data in Figure \ref{case10}. The structure of the
modules $\Delta_{10}(521^3)$, $\Delta_{10}(4321)$ and
$\Delta_{10}(4211)$, follows exactly as in the preceding cases for
partitions of $n$ and $n-2$. For $\Delta_{10}(332)$ we also need to
note that $(332)\not\subset(521^3)$, and so
$L_{10}(521^3)$ cannot occur.  

To see that $L_{10}(521^3)$ cannot occur anywhere else it is enough to
note (by Proposition \ref{projelim}) that $\Delta_9(421^3)$ is
projective. Similarly $L_{10}(4321)$ cannot occur in the standards
$\Delta_{10}(22)$ and $\Delta_{10}(0)$ as $\Delta_8(431)$ is
projective. The structure of $\Delta_{10}(0)$ and $\Delta_{10}(22)$
then follows by localisation to $n=8$.

The only remaining module is $\Delta_{10}(321)$. It is clear that this
must have at most the four factors shown. The multiplicities of
$L_{10}(4211)$ and $L_{10}(332)$ must be $1$ by localisation to the case
$n=8$. It remains to show that the final factor has multiplicity $1$,
that there is a map to the module from $\Delta_{10}(4321)$, 
and that the module structure is as shown.

Consider $\res_{10}L_{10}(4321)=\res_{10}\Delta_{10}(4321)$.  By
(\ref{resrule}), and the simplicity of standard modules
$\Delta_n(\lambda)$ when $\lambda\vdash n$, this has simple factors
\begin{equation}
\label{10simfac}
L_9(432)\quad\quad L_9(4311)\quad\quad L_9(4221)\quad\quad
L_9(3321).
\end{equation}
If we consider $\res_{10}\Delta_{10}(4211)$ and
$\res_{10}\Delta_{10}(332)$, using the structure of
$\Delta_{10}(4211)$ and $\Delta_{10}(332)$ given above, it is easy to
show that neither $\res_{10}L_{10}(4211)$ nor $\res_{10}L_{10}(332)$
contain any of the factors in (\ref{10simfac}).  Comparing with
(\ref{9res}), (\ref{9lower}), and (\ref{9upper}) we see that
$$[\Delta_{10}(321):L(4321)]\leq 1.$$
Further, either this simple does occur, or the simples in
(\ref{10simfac}) all occur in $\res_{10}L_{10}(321)$.

By (\ref{ressimple}) we have that $\res_{10}L_{10}(4211)$ contains
$L_9(421)$ and $\res_{10}L_{10}(332)$ contains $L_9(331)$. (In fact,
all factors can be easily determined.)  Comparing with (\ref{9upper})
we see that both $L_{10}(4211)$ and $L_{10}(331)$ must occur above the
simple whose restrictions contribute the terms in (\ref{10simfac}).
Thus $L_{10}(4321)$ must occur and forms the socle of
$\Delta_{10}(321)$. This completes our verification of the various
claims above, and so the block structure is as shown in Figure
\ref{case10}.

\subsection{The case $n=11$}

\begin{figure}
\includegraphics{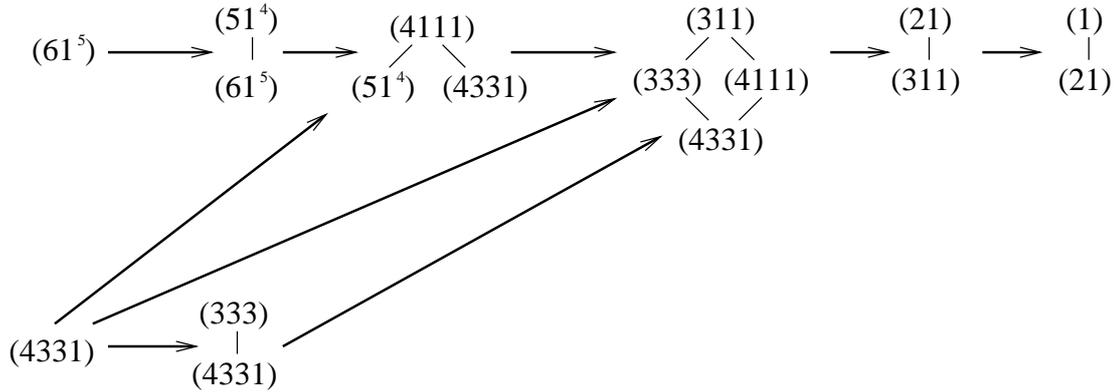}
\caption{The block containing $(1)$ when $n=11$}\label{case11}
\end{figure}

We claim that the structure of the block containing $(1)$ is given by
the data in Figure \ref{case11}. The structure of the modules
$\Delta_{11}(61^5)$, $\Delta_{11}(4331)$, $\Delta_{11}(51^4)$ and
$\Delta_{11}(333)$ follows as for $n=10$ for partitions of $n$ and
$n-2$.

The modules $\Delta_{10}(51^5)$ and $\Delta_8(431)$ are
projective. Therefore by Proposition \ref{projelim} the simple
$L_{11}(61^5)$ cannot occur in any of the remaining standards, and
$L_{11}(4331)$ cannot occur in $\Delta_{11}(21)$ or
$\Delta_{11}(1)$. The structure of these latter two modules now
follows by localisation to the case $n=9$. To see that the structure
of $\Delta_{11}(41^3)$ is as illustrated follows from
(\ref{squarestuff}) and localisation.

The only remaining module is $\Delta_{11}(311)$. By localising this
contains $L_{11}(31^2)$, $L_{11}(3^3)$ and $L_{11}(41^3)$, all with
multiplicity. We have eliminated all other possible factors except
$L_{11}(4331)$. We proceed as for the module $\Delta_{10}(321)$
above. Restriction of $\Delta_{11}(311)$ contains $L_{10}(433)$ with
multiplicity one, and this can only arise in the restriction of
$L_{11}(4311)$. Arguing as in the $n=10$ case we also see that
$L_{11}(4311)$ must coincide with the socle of $\Delta_{11}(311)$, and
so we are done.

\subsection{The case $n=12$}

\begin{figure}
\center{\includegraphics[angle=90]{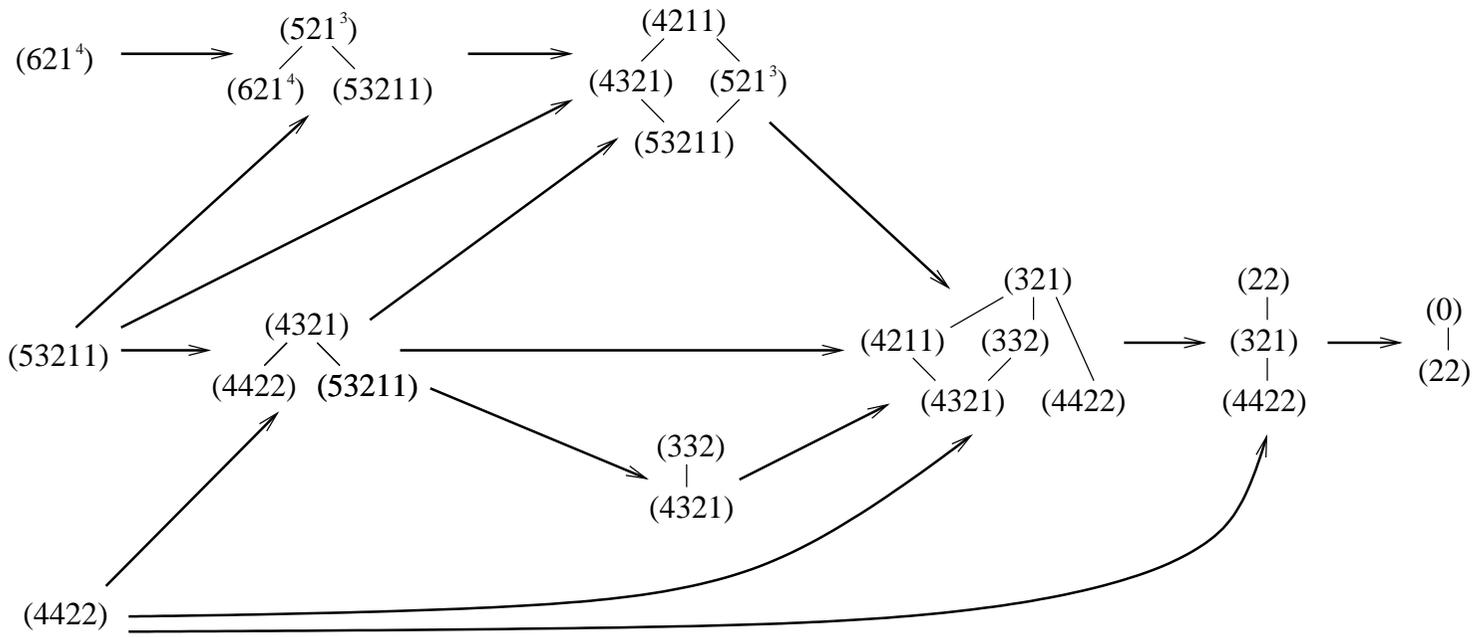}}
\caption{The block containing $0$ when $n=12$}\label{case12}
\end{figure}

We claim that the structure of the block containing $(0)$ is given by
the data in Figure \ref{case12}.
As usual, the structure of the modules labelled by partitions of $n$
and $n-2$ is straightforward. The simple $L_{12}(621^4)$ cannot occur
anywhere else as $\Delta_{11}(521^4)$ is projective.

The modules $\Delta_{10}(42211)$ and $\Delta_8(2^4)$ are both
projective. Therefore $L_{12}(53211)$ cannot occur in any standard
labelled by a partition of $6$ or smaller, while $L_{12}(4422)$ cannot
occur in $\Delta_{12}(0)$. The structure of $\Delta_{12}(0)$ is then
clear by localisation.

Next consider $\Delta_{12}(332)$. This cannot contain $L_{12}(53211)$
as $(332)\not\subset(53211)$, so it is enough by localisation to
verify that $L_{12}(4422)$ cannot occur. But this is clear, as the
restriction of this simple contains $L_{11}(4421)$, which is not in
the same block as any of the standard modules in the restriction of
$\Delta_{12}(332)$. 

It remains to determine the structure of $\Delta_{12}(321)$ and
$\Delta_{12}(22)$. In each case we know the multiplicity of all
composition factors by localisation and the remarks above, except for
$L_{12}(4422)$. Using Proposition \ref{largerhom} we see that we have
a non-zero homomorphism from $\Delta_{12}(4422)$ into each of the two
standards, and so the multiplicity of $L_{12}(4422)$ in each case is
at least one. The restriction of each of these standards contains
precisely one copy of $L_{11}(4322)$ and $L_{11}(4421)$, and so
$L_{12}(4422)$ must occur with multiplicity one in each standard. 

To determine the location of $L_{12}(4422)$ in $\Delta_{12}(321)$ note
that it cannot occur below any composition factor other than
$L_{12}(321)$, as this would contradict the existence of homomorphisms
from $\Delta_{12}(4211)$ and $\Delta_{12}(332)$ into
$\Delta_{12}(321)$. Therefore the structure of this module must be as
shown. For $\Delta_{12}(22)$, we have that the simple $L_{11}(32)$ in
the restriction of $L_{12}(321)$ occurs above $L_{12}(4421)$, and
hence the structure of $\Delta_{12}(22)$ must be as shown.

\subsection{A comparison with Kazhdan-Lusztig polynomials}

Suppose that $W$ is a Weyl group, with associated affine Weyl group
$W_p$.  Soergel has shown \cite{soergel1,soergel2} that (provided $p$
is not too small) the value of the parabolic Kazhdan-Lusztig
polynomials $n_{\lambda\mu}$ (evaluated at $v=1$) associated to
$W\subset W_p$ determine the multiplicity of the standard module
$\Delta_q(\lambda)$ in the indecomposable $T_q(\mu)$ for the quantum
group $U_q$ associated to $W$ where $q$ is a $p$th root of unity.

In the case of the quantum general linear group, Ringel duality
\cite{erdqh} translates this into a result about decomposition numbers
for the Hecke algebra of type $A$, where now $\mu$ labels a simple
module and $\lambda$ a Specht module. Further, Rouquier has
conjectured that the coefficient of $v^t$ occurring in $n_{\lambda\mu}$
should correspond to the multiplicity of the simple $D^{\mu}$ in the
$t$th layer of the Jantzen filtration of the Specht module
$S^{\lambda}$.

In this spirit, we can compare our results in this section with the
polynomials in Figure \ref{klcalc} for $n\leq 12$. We see that in each
case, the value of $n_{\lambda,\mu}(1)$ from Figure \ref{klcalc} is
exactly the multiplicity of $L_n(\mu)$ in $\Delta_n(\lambda)$, and
that there is a filtration of $\Delta_{n}(\lambda)$ corresponding to
the powers of $v$ occurring in the polynomials for the
$L_n(\mu)$s. This, together with the other Lie-like phenomena we have
observed leads us to ask

\begin{qn} (i) For the Brauer algebra with Kazhdan-Lusztig polynomials
  as defined in Section \ref{klpoly}, is it true for weights in an
  alcove that
$$[\Delta_n(\lambda):L_n(\mu)]=n_{\lambda,\mu}(1)\ ?$$ (ii) Is there a
  (Jantzen?) filtration of $\Delta_n(\lambda)$ such that the
  multiplicity of a simple $L_n(\mu)$ in the $t$th layer is given by
  the coefficient of $v^t$ in $n_{\lambda,\mu}$?
\end{qn}

As we have noted, the results in this section answer both parts in the
affirmative when $n\leq 12$ and $\delta=1$.

\providecommand{\bysame}{\leavevmode\hbox to3em{\hrulefill}\thinspace}
\providecommand{\MR}{\relax\ifhmode\unskip\space\fi MR }
\providecommand{\MRhref}[2]{%
  \href{http://www.ams.org/mathscinet-getitem?mr=#1}{#2}
}
\providecommand{\href}[2]{#2}

\end{document}